# NEGATIVE SOBOLEV REGULARITY AND OPTIMAL DECAY FOR THE FRACTIONAL NAVIER–STOKES–POISSON SYSTEM

YUAN XU AND JUNHAO ZHANG

ABSTRACT. We establish global well-posedness and sharp higher-order decay estimates for the fractional compressible Navier–Stokes–Poisson system in $\mathbb{R}^3$ with fractional Lamé–Navier dissipation $\mathfrak{L}^\alpha u$ with $\alpha \in [1/2, 1]$, and Riesz coupling $\Lambda^{2\beta}\phi = \rho - \bar{\rho}$ with $\beta \in [0, 1]$, where $\Lambda := (-\Delta)^{1/2}$ and $\bar{\rho} > 0$ denotes the constant background density. For small perturbations of a constant equilibrium in $H^{N\alpha}$ with $N \geq 3$ and $N\alpha > 5/2$, we prove the global existence and uniqueness of classical solutions. In addition, we propagate the negative Sobolev norms with uniform-in-time bounds and derive the higher-order decay estimates, provided that the initial perturbation additionally belongs to $\dot{H}^{-\vartheta}$ for some $\vartheta \in [0, 3/2)$. Moreover, we establish the optimality of these decay rates within the prescribed initial data class for the velocity and the Riesz potential in their natural energy norms. The density decay exponents are also optimal whenever the available regularity permits the full Riesz-induced gain. The analysis relies on a nonlinear energy structure adapted to both fractional dissipation and Riesz interaction, together with Sobolev interpolation. The optimality of the decay is proved by constructing a special family of smooth initial data localized at low frequencies and proving matching upper and lower decay bounds for the resulting solutions. These solutions demonstrate that the decay dictated by the low-frequency linear dynamics persists under nonlinear evolution, including the additional density decay induced by the Riesz coupling.

## 1. Introduction

The evolution of a compressible fluid in $\mathbb{R}^3$ with fractional Lamé–Navier dissipation and Riesz interaction is described by the following Cauchy problem:

$$\begin{cases} \partial_t \rho + \operatorname{div}(\rho u) = 0, \\ \partial_t(\rho u) + \operatorname{div}(\rho u \otimes u) + \nabla p(\rho) + \rho \nabla \phi = -\mathfrak{L}^\alpha u, \\ \Lambda^{2\beta}\phi = \rho - \bar{\rho}, \\ (\rho, u)|_{t=0} = (\rho_0, u_0). \end{cases} \tag{1.1}$$

Here $\rho = \rho(t, x) \geq 0$ and $u = u(t, x)$ denote the density and velocity of the fluid at time $t \geq 0$ and position $x \in \mathbb{R}^3$, respectively. The function $\phi = \phi(t, x)$ is the self-consistent electric potential generated by the density fluctuation. The pressure $p$ is given by

$$p(\rho) = K\rho^\gamma,$$

where $K > 0$ is a pressure constant and $\gamma > 1$ is the adiabatic exponent. We set $\Lambda := (-\Delta)^{1/2}$. The positive fractional Lamé–Navier operator is defined by

$$\widehat{\mathfrak{L}^\alpha u}(\xi) := |\xi|^{2\alpha}\big(\nu_{\mathrm{T}}\mathbb{P}_{\mathrm{T}}(\xi) + \nu_{\mathrm{L}}\mathbb{P}_{\mathrm{L}}(\xi)\big)\widehat{u}(\xi) := \mathbb{B}_\alpha(\xi)|\xi|^{2\alpha}\widehat{u}(\xi), \tag{1.2}$$

where $\nu_{\mathrm{T}} = \mu^\alpha > 0$, $\nu_{\mathrm{L}} = (2\mu + \lambda)^\alpha > 0$, and $\mathbb{B}_\alpha = \nu_{\mathrm{T}}\mathbb{P}_{\mathrm{T}} + \nu_{\mathrm{L}}\mathbb{P}_{\mathrm{L}}$. The longitudinal and transverse Helmholtz projections are defined as

$$\mathbb{P}_{\mathrm{L}} v := \nabla\Delta^{-1} \operatorname{div} v, \qquad \mathbb{P}_{\mathrm{T}} v := (\mathbb{I} - \mathbb{P}_{\mathrm{L}})v,$$

where $\mathbb{I}$ denotes the identity operator on the relevant function space. The parameters $\alpha \in [1/2, 1]$ and $\beta \in [0, 1]$ specify the fractional powers in the dissipation and the Riesz interaction, respectively. The constant $\bar{\rho} > 0$ is the uniform background density. We refer to $\Lambda^{\beta}\phi$ as the associated Riesz potential, with the endpoint $\beta = 0$ understood as $\rho - \bar{\rho}$. For simplicity, we take $\bar{\rho} = 1$ and $K\gamma = 1$, so that $p'(1) = 1$.

The system (1.1) involves both fractional Lamé–Navier dissipation and Riesz interaction. We first recall that the viscous operator in (1.2) reduces to the classical positive Lamé–Navier operator $\mathcal{L}u = -\mu\Delta u - (\mu + \lambda)\nabla \operatorname{div} u$ when $\alpha = 1$. Physically, $-\mathcal{L}u$ represents the divergence of the viscous stress tensor in compressible Newtonian fluids and the elastic restoring force in homogeneous isotropic solids, with its transverse and longitudinal components corresponding to shear and compressional modes, respectively. For the fractional Lamé–Navier operator in (1.2), we adopt Scott's positive spectral realization [28]. The analytic properties used below are summarized in Appendix A. Such fractional extensions also arise naturally in the modeling of spatially nonlocal materials. Nonlocal continuum theories allow the constitutive response at a point to depend on deformation away from that point (see Eringen [11]). Within this framework, Evgrafov and Bellido investigated the connection between Eringen's nonlocal elasticity and fractional elasticity [12]. Šilhavý developed an isotropic fractional model preserving the transverse–longitudinal decomposition [29]. These developments support the interpretation of fractional Lamé–Navier operators as continuum models of power-law spatial nonlocality.

In addition to fractional viscous dissipation, the system (1.1) includes a Riesz coupling that, for $0 < \beta \leq 1$, models repulsive nonlocal interactions generated by the density fluctuation $\rho - \bar{\rho}$ relative to a uniform background. At $\beta = 1$, the Riesz interaction reduces to the classical Coulomb law, with a potential kernel proportional to $|x|^{-1}$ and an inverse-square repulsive force. For $0 < \beta < 1$, it provides a power-law generalization of Coulomb repulsion, with a potential kernel that is more singular at the origin and decays faster at infinity than the Coulomb kernel. Riesz potentials arise in the statistical mechanics of systems with long-range interactions [6], including Riesz gases [19]. The particle, kinetic, and hydrodynamic descriptions of such collective behavior are reviewed by Carrillo, Choi, and Perez [7].

We now turn to related mathematical results for fluid and kinetic equations involving nonlocal effects. Global existence has been established for incompressible models with fractional dissipation, such as the rapidly rotating Navier–Stokes system [1], the Boussinesq system [31], and the generalized SQG equations [17]. Time-decay estimates have also been obtained for the Oldroyd-B model with fractional stress diffusion [34] and for the fractional Fokker–Planck model [3]. For compressible Navier–Stokes equations with fractional viscous dissipation, Wang and Zhang [35] studied the Cauchy problem with $\alpha \in (1/2, 1]$. They established the nonlinear stability of the constant equilibrium state and proved the global existence of small classical solutions within an $H^4$ framework. Furthermore, under an additional $L^1$ assumption on the initial perturbation, they obtained an $L^2$ upper bound of order $(1+t)^{-3/(4\alpha)}$ for the density-velocity perturbation. In a subsequent study [36], they refined this theory by relaxing the regularity requirement from $H^4$ to the fractional Sobolev space $H^{3\alpha+1}$.

The compressible Navier–Stokes system with classical viscosity and Riesz interaction has been studied in both the repulsive and attractive regimes. In the repulsive regime with $\beta \in (0, 1)$, Świerczewska-Gwiazda, Xu, and Zhang [33] established global smooth solutions near a constant equilibrium, $L^2$ and $L^\infty$ decay estimates, and a quantitative global-in-time inviscid limit under suitable smallness assumptions. In the attractive regime with $\beta \in (3/8, 1/2)$, Carrillo, Duan, Wróblewska-Kamińska, and Zhang [8] proved global existence of strong solutions and asymptotic stability of compactly supported steady states for the one-dimensional problem with physical vacuum boundaries under suitable small perturbations.

For the compressible Navier–Stokes–Poisson system with both fractional viscous dissipation and Riesz interaction, Wang and Zhang [37] considered $\alpha \in (1/2, 1]$ and $\beta \in [1/2, 1]$. Based on a detailed spectral analysis of the linearized semigroup and time-weighted energy estimates, they proved the $\dot{H}^s$ upper bounds of orders $(1+t)^{-(3+2s)/(4\alpha)}$ and $(1+t)^{-(3+2s-2\beta)/(4\alpha)}$ for $\rho-1$ and $(u, \Lambda^\beta\phi)$, respectively. These estimates show that the fractional electric coupling slows down the decay of $(u, \Lambda^\beta\phi)$, compared with the fractional compressible Navier–Stokes system. However, that work does not cover the weak-coupling regime $\beta < 1/2$.

In this paper, we study global well-posedness and optimal time decay for the system (1.1). We establish the global well-posedness of small perturbations of a constant equilibrium for the compressible Navier–Stokes system with the fractional Lamé–Navier operator in (1.2) and repulsive Riesz interaction, for $\alpha \in [1/2, 1]$ and $\beta \in [0, 1]$. Under an additional finite negative Sobolev assumption, we derive a hierarchy of higher-order time-decay estimates, including an additional density decay gain of $\beta/(2\alpha)$ whenever the available regularity permits. We further establish the optimality of the decay exponents within the prescribed initial-data class by constructing smooth initial data localized at low frequencies and proving matching upper and lower decay bounds for the resulting nonlinear solutions. This optimality holds for the velocity and the Riesz potential, and for the density whenever its decay is not limited by the available regularity. These results are summarized in the following theorem.

**Theorem 1.1.** *Let $\alpha \in [1/2, 1]$, $\beta \in [0, 1]$, and $\mu > 0$, with $2\mu + \lambda > 0$. Let $N \geq 3$ be an integer satisfying $N\alpha > 5/2$.*
(i) *There exists a constant $\delta_0 > 0$ such that if the initial data satisfy*

$$\sum_{k=0}^{N} \left( \|\Lambda^{k\alpha}(\rho_0 - 1)\|_{L^2}^2 + \|\Lambda^{k\alpha}u_0\|_{L^2}^2 + \|\Lambda^{k\alpha+\beta}\phi_0\|_{L^2}^2 \right) \leq \delta_0^2, \tag{1.3}$$

*then the system* (1.1) *admits a unique global solution $(\rho, u, \Lambda^\beta\phi)$ satisfying, for all $t \geq 0$,*

$$\begin{aligned} &\sum_{k=0}^{N} \left( \|\Lambda^{k\alpha}(\rho - 1)(t)\|_{L^2}^2 + \|\Lambda^{k\alpha}u(t)\|_{L^2}^2 + \|\Lambda^{k\alpha+\beta}\phi(t)\|_{L^2}^2 \right) \\ &\quad + \int_0^t \sum_{k=1}^{N} \left( \|\Lambda^{k\alpha}(\rho - 1)(\tau)\|_{L^2}^2 + \|\Lambda^{k\alpha+\beta}\phi(\tau)\|_{L^2}^2 \right) + \sum_{k=1}^{N+1} \|\Lambda^{k\alpha}u(\tau)\|_{L^2}^2 \,\mathrm{d}\tau \\ &\leq C \sum_{k=0}^{N} \left( \|\Lambda^{k\alpha}(\rho_0 - 1)\|_{L^2}^2 + \|\Lambda^{k\alpha}u_0\|_{L^2}^2 + \|\Lambda^{k\alpha+\beta}\phi_0\|_{L^2}^2 \right). \end{aligned} \tag{1.4}$$

(ii) *If, in addition,*

$$M_0 := \|\Lambda^{-\vartheta}(\rho_0 - 1)\|_{L^2} + \|\Lambda^{-\vartheta}u_0\|_{L^2} + \|\Lambda^{-\vartheta+\beta}\phi_0\|_{L^2} < +\infty, \tag{1.5}$$

*for some $\vartheta \in [0, 3/2)$, we have, for any $t \geq 0$,*

$$\|\Lambda^{-\vartheta}(\rho - 1)(t)\|_{L^2} + \|\Lambda^{-\vartheta}u(t)\|_{L^2} + \|\Lambda^{-\vartheta+\beta}\phi(t)\|_{L^2} \leq C\widetilde{M}_0, \tag{1.6}$$

*and, for $l = 0, \ldots, N-1$,*

$$\sum_{k=l}^{N} \left( \|\Lambda^{k\alpha}(\rho - 1)(t)\|_{L^2} + \|\Lambda^{k\alpha}u(t)\|_{L^2} + \|\Lambda^{k\alpha+\beta}\phi(t)\|_{L^2} \right) \leq C\widetilde{M}_0(1+t)^{-\left(\frac{\vartheta}{2\alpha} + \frac{l}{2}\right)}, \tag{1.7}$$

*while, for $l = 0, \ldots, N-2$,*

$$\|\Lambda^{l\alpha}(\rho - 1)(t)\|_{L^2} \leq C\widetilde{M}_0(1+t)^{-\frac{\vartheta}{2\alpha} - \frac{l}{2} - \frac{\min\{\beta, (N-1-l)\alpha\}}{2\alpha}}. \tag{1.8}$$

*In the above estimates, $\widetilde{M}_0$ is defined by*

$$\widetilde{M}_0 := \begin{cases} M_0 + \delta_0, & 0 \le \vartheta \le 3/2 - \alpha, \\ (M_0+\delta_0)(1+M_0+\delta_0), & 3/2-\alpha < \vartheta < 3/2. \end{cases}$$

(iii) *For each fixed $\vartheta \in [0, 3/2)$, the decay exponent in* (1.7) *is optimal for $\|\Lambda^{l\alpha} u(t)\|_{L^2}$ and $\|\Lambda^{l\alpha+\beta}\phi(t)\|_{L^2}$, for $l = 0, \dots, N-1$. The density decay exponent in* (1.8) *is also optimal for $l = 0, \dots, N-2$ with $(N-1-l)\alpha \ge \beta$.*

**Remark 1.2.** *Within the common parameter and derivative ranges, Theorem 1.1 yields faster decay rates with an additional decay factor of $(1+t)^{-(\vartheta+\beta-3/2)/(2\alpha)}$ than the rates established in [37] when $\vartheta \in (3/2-\beta, 3/2)$. For the density, this comparison requires $(N-1-l)\alpha \ge \beta$ at derivative order $l\alpha$.*

**Remark 1.3.** *The parameter range includes the weak-coupling regime $\beta \in [0, 1/2)$, and the model recovers two familiar coupling regimes at the endpoints $\beta = 0$ and $\beta = 1$. When $\beta = 1$, the interaction law becomes the classical Poisson relation $-\Delta\phi = \rho - 1$. When $\beta = 0$, the relation $\phi = \rho - 1$ yields $\rho\nabla\phi = \nabla(\rho^2/2)$. Consequently, the coupling term can be absorbed into the pressure. Therefore, the present framework treats both the classical nonlocal Poisson coupling and the local pressure-type coupling in a unified manner.*

**Remark 1.4.** *For the initial data, the condition imposed on the Riesz potential is equivalent to an additional low-frequency assumption on the density perturbation. This additional control is precisely what yields the enhanced density decay in* (1.8)*, compared with the fractional compressible Navier–Stokes system without electric coupling. More precisely, strictly below the highest regularity threshold, the Riesz relation $\Lambda^{2\beta}\phi = \rho - 1$ introduces an extra decay factor of $(1+t)^{-\beta/(2\alpha)}$. However, at the highest density derivative order $l = N-2$, the full gain is constrained by the available regularity when $\beta > \alpha$. This restriction accounts for the appearance of the term $\min\{\beta, (N-1-l)\alpha\}$ in* (1.8)*. In the classical case $\alpha = \beta = 1$, this recovers the well-known additional one-half decay rate for the density in the Navier–Stokes–Poisson system. In contrast, when $\beta = 0$, the additional gain vanishes, which is consistent with the fractional compressible Navier–Stokes scaling.*

**Remark 1.5.** *The optimality asserted in Theorem 1.1 is understood within the negative Sobolev framework. More precisely, for fixed $N$ and $\vartheta$, increasing any decay exponent in Theorem 1.1*(iii) *by any positive amount makes the corresponding estimate fail for some initial data satisfying conditions* (1.3) *and* (1.5)*. In the classical case $\alpha = \beta = 1$, these exponents agree with those in [39] in the common derivative ranges. The present result complements those estimates by establishing the optimality of the corresponding exponents within the initial-data class considered here.*

To prove the optimality in Theorem 1.1(iii), we construct a family of special solutions with matching upper and lower decay estimates. The auxiliary regularity index $N_*$ below is used only to construct a smooth special solution at the prescribed derivative order $\sigma$.

**Theorem 1.6.** *Let $\alpha \in [1/2, 1]$, $\beta \in [0, 1]$, $0 \le \vartheta < 3/2$, $\sigma \ge 0$, and $0 < \varepsilon \le 3/2 - \vartheta$. Let $N_*$ be an integer satisfying*

$$(N_* - 1)\alpha > \max\{\sigma + \beta + 2,\ \beta + 3,\ 2\alpha + 2\}. \tag{1.9}$$

*Then there exists $\delta_* > 0$ such that, for any $0 < \delta \le \delta_*$, one can choose a class of real-valued initial data $\rho_0^{\delta,\varepsilon} = 1 + \varrho_0^{\delta,\varepsilon}$, $u_0^{\delta,\varepsilon}$, $\Lambda^\beta \phi_0^{\delta,\varepsilon}$ satisfying* (1.3) *and* (1.5) *with $N = N_*$ and negative Sobolev index $\vartheta$. For these initial data, Theorem 1.1*(i) *ensures the existence of a unique global solution*

$(\rho^{\delta,\varepsilon}, u^{\delta,\varepsilon}, \Lambda^\beta \phi^{\delta,\varepsilon})$ *to* (1.1)*. Moreover, for any* $t \geq 0$*, the corresponding solution satisfies*

$$c\delta(1+t)^{-\frac{\sigma+\beta+\vartheta+\varepsilon}{2\alpha}} \leq \|\Lambda^\sigma(\rho^{\delta,\varepsilon}-1)(t)\|_{L^2} \leq C\delta(1+t)^{-\frac{\sigma+\beta+\vartheta+\varepsilon}{2\alpha}}, \tag{1.10}$$

$$c\delta(1+t)^{-\frac{\sigma+\vartheta+\varepsilon}{2\alpha}} \leq \|\Lambda^\sigma u^{\delta,\varepsilon}(t)\|_{L^2} \leq C\delta(1+t)^{-\frac{\sigma+\vartheta+\varepsilon}{2\alpha}}, \tag{1.11}$$

$$c\delta(1+t)^{-\frac{\sigma+\vartheta+\varepsilon}{2\alpha}} \leq \|\Lambda^{\sigma+\beta}\phi^{\delta,\varepsilon}(t)\|_{L^2} \leq C\delta(1+t)^{-\frac{\sigma+\vartheta+\varepsilon}{2\alpha}}. \tag{1.12}$$

**Remark 1.7.** *The initial perturbations constructed in the preceding theorem belong to* $H^m$ *for every* $m \geq 0$*. Nevertheless, the low-frequency profiles of* $u_0$ *and* $\Lambda^\beta\phi_0$ *approach the integrability threshold associated with* $\dot{H}^{-\vartheta}$ *as* $\varepsilon \to 0$*. At the derivative orders specified in Theorem 1.1*(iii)*, the decay exponents of these solutions exceed the corresponding exponents by* $\varepsilon/(2\alpha)$*. Since* $\varepsilon > 0$ *can be arbitrarily small, the lower bounds exclude any positive increase in those exponents throughout the admissible initial-data class. Thus, arbitrarily high positive Sobolev regularity alone does not improve these decay exponents.*

By virtue of the Hardy–Littlewood–Sobolev inequality (Lemma A.5), we have the continuous embedding $L^p(\mathbb{R}^3) \subset \dot{H}^{-\vartheta}(\mathbb{R}^3)$ for $p \in (1,2]$ with $\vartheta = 3\,(1/p - 1/2) \in [0, 3/2)$. Consequently, the decay estimates in Theorem 1.1(ii) yield the following $L^p$–$L^2$ upper bounds.

**Corollary 1.8.** *Suppose that the hypotheses of Theorem 1.1*(i) *hold. If, instead of the negative Sobolev assumption in part* (ii)*, we assume*

$$M_1 := \|(\rho_0 - 1)\|_{L^p} + \|u_0\|_{L^p} + \|\Lambda^\beta\phi_0\|_{L^p} < +\infty,$$

*for some* $p \in (1,2]$*, then, for all* $t \geq 0$ *and* $l = 0, \ldots, N-1$*,*

$$\sum_{k=l}^{N}\Big(\|\Lambda^{k\alpha}(\rho-1)(t)\|_{L^2} + \|\Lambda^{k\alpha}u(t)\|_{L^2} + \|\Lambda^{k\alpha+\beta}\phi(t)\|_{L^2}\Big) \leq C(M_1,\delta_0)(1+t)^{-\frac{3}{2\alpha}\left(\frac{1}{p}-\frac{1}{2}\right)-\frac{l}{2}}.$$

*Moreover, for* $l = 0, \ldots, N-2$*,*

$$\|\Lambda^{l\alpha}(\rho-1)(t)\|_{L^2} \leq C(M_1,\delta_0)(1+t)^{-\frac{3}{2\alpha}\left(\frac{1}{p}-\frac{1}{2}\right)-\frac{l}{2}-\frac{\min\{\beta,(N-1-l)\alpha\}}{2\alpha}}.$$

To place our decay and optimality results in the context of the classical theory, we recall several developments concerning the global existence and large-time behavior of compressible viscous flows. For the classical compressible Navier–Stokes equations in the whole space, Matsumura and Nishida established global existence and asymptotic stability of small perturbations of constant equilibria through energy estimates [25], and derived algebraic decay using linearized spectral analysis and nonlinear energy estimates [24]. Matsumura [23] subsequently established global existence and decay through time-weighted energy estimates and Sobolev interpolation. Through diffusion-wave analysis, Hoff and Zumbrun [15] established matching $L^2$ upper and lower decay bounds for the joint density-momentum perturbation arising from suitable small, localized initial data. For further results on diffusion waves and pointwise decay, see [16, 22]. Guo and Wang [13] developed a pure energy method that combines the propagation of negative Sobolev norms with interpolation among energy estimates to derive decay rates for small solutions and their higher-order spatial derivatives. More recently, Brandolese, Shou, Xu, and Zhang [5] established a sharp decay characterization for small solutions in a critical Besov framework, relating matching upper and lower bounds to the low-frequency regularity and nondegeneracy of the initial perturbation. This characterization builds on Brandolese's earlier work [4] on low-frequency criteria for sharp decay of solutions to the incompressible Navier–Stokes equations. The estimates in [5], combined with suitable low-frequency initial data and Sobolev interpolation, also imply the optimality of the joint Sobolev decay exponents obtained by Guo and Wang [13] for the classical compressible Navier–Stokes equations.

A more direct classical counterpart of the present model is the compressible Navier–Stokes–Poisson system, which corresponds to the case $\alpha = \beta = 1$ in (1.1). For this system, Li, Matsumura, and Zhang [21] established global classical solutions for sufficiently small initial data in $H^N \cap L^1$ with $N \geq 4$. They obtained $L^2$ upper bounds of order $(1+t)^{-3/4}$ for the density perturbation and $(1+t)^{-1/4}$ for the momentum and electric field. In addition, they proved matching $L^2$ upper and lower decay bounds under the additional condition that $|\widehat{\rho_0 - \bar{\rho}}(\xi)| \geq c_0 > 0$ for all sufficiently small $|\xi|$. Other global existence results can be found in [10, 14, 38]. The natural energy of the system contains an electric-field contribution, so finite-energy perturbations require $\nabla\phi_0 \in L^2$, which further implies $\Lambda^{-1}(\rho_0 - 1) \in L^2$ by the Poisson equation. Motivated by this energy structure, Wang [39] incorporated the electric field into the assumptions on the initial perturbation. He considered initial data in $H^N \cap \dot{H}^{-s}$ with $N \geq 3$ and $s \in [0, 3/2)$, together with the neutrality condition. Assuming that the $H^3$ norm of this initial perturbation is sufficiently small, he established a unique global solution and obtained decay estimates through a refined pure energy method combining the propagation of negative Sobolev norms with interpolation among energy estimates. In particular, the $L^2$ upper bounds are of order $(1+t)^{-(s+1)/2}$ for the density perturbation and $(1+t)^{-s/2}$ for the velocity and electric field, with corresponding estimates for higher-order spatial derivatives. When the initial perturbation belongs to $L^p$, with $p \in (1, 2]$, Wang established $L^2$ upper bounds of order $(1+t)^{-3/2(1/p-1/2)-1/2}$ for the density perturbation and $(1+t)^{-3/2(1/p-1/2)}$ for the velocity and electric field. These estimates reflect the additional low-frequency control of the density provided by the electric-field assumption.

We next outline the proof strategy. The standard energy estimate for (1.1) gives direct dissipation only for the velocity through the fractional Lamé–Navier dissipation $\mathfrak{L}^\alpha u$. The density and the potential are not directly dissipated at this level. Moreover, a direct estimate of the force $\nabla\phi$ in terms of the energy variable $\Lambda^\beta\phi$ involves $1-\beta$ additional derivatives. This derivative gap increases as $\beta$ decreases, requiring particular care in the treatment of the potential coupling. To overcome these difficulties, we first introduce the sound-speed variable

$$\mathfrak{c} = \frac{2}{\gamma - 1}\left((1+\varrho)^{\frac{\gamma-1}{2}} - 1\right), \qquad \varrho = \rho - 1,$$

which symmetrizes the hyperbolic part and avoids derivative loss from the pressure term. Under the smallness assumption, the Sobolev norms of $\mathfrak{c}$ and $\varrho$ are equivalent. The key new ingredient is the construction of an energy functional adapted to the fractional order. More precisely, besides the standard energy of $(\mathfrak{c}, u, \Lambda^\beta\phi)$, we introduce the interactive term

$$\int_{\mathbb{R}^3} \Lambda^{k\alpha} u \cdot \nabla \Lambda^{k\alpha + 2\alpha - 2} \varrho \, \mathrm{d}x.$$

This term is designed to match the differential order $2\alpha$ of the viscous dissipation and therefore recovers the hidden dissipation of the density at the level $\|\Lambda^{(k+1)\alpha}\varrho\|_{L^2}^2$. By using the Riesz relation $\Lambda^{2\beta}\phi = \varrho$, the same estimate also yields the dissipation of the potential in the form $\|\Lambda^{(k+1)\alpha+\beta}\phi\|_{L^2}^2$. Thus, the dissipation of $\varrho$ and $\Lambda^\beta\phi$ is recovered directly at the nonlinear energy level.

At the highest order, the treatment of the electric term is especially delicate. We use the continuity equation $\rho\,\mathrm{div} u = -\partial_t\rho - u \cdot \nabla\rho$, together with the Riesz relation, to extract a time derivative of the weighted potential energy

$$\left\| \frac{1}{\sqrt{\rho}} \Lambda^{N\alpha+\beta}\phi \right\|_{L^2}^2$$

from the coupling term $(\Lambda^{N\alpha}\nabla\phi, \Lambda^{N\alpha} u)$. The remaining commutator and transport terms are controlled by the available energy and dissipation, without requiring additional regularity. Since $\rho$ remains close to 1, this weighted energy is equivalent to the usual highest-order potential energy. This yields the highest-order estimate without derivative loss. Combining this estimate with the lower-order energy estimates and sufficiently small multiples of the interactive dissipation estimates, we construct a modified energy $\mathcal{E}_N^l$ for each $l = 0, \ldots, N-1$. Under the smallness assumption, there exist constants $c$, $C > 0$, independent of $t$, such that

$$c\mathcal{E}_N^l(t) \leq \sum_{k=l}^{N} \left( \|\Lambda^{k\alpha}\varrho(t)\|_{L^2}^2 + \|\Lambda^{k\alpha}u(t)\|_{L^2}^2 + \|\Lambda^{k\alpha+\beta}\phi(t)\|_{L^2}^2 \right) \leq C\mathcal{E}_N^l(t).$$

The precise construction of $\mathcal{E}_N^l$ is given in Section 4. The nonlinear error terms can then be absorbed by the dissipation, giving a Lyapunov-type inequality

$$\frac{\mathrm{d}}{\mathrm{d}t}\mathcal{E}_N^l(t) + c\mathcal{D}_N^l(t) \leq 0,$$

where

$$\mathcal{D}_N^l(t) := \sum_{k=l+1}^{N} \left( \|\Lambda^{k\alpha}\varrho(t)\|_{L^2}^2 + \|\Lambda^{k\alpha+\beta}\phi(t)\|_{L^2}^2 \right) + \sum_{k=l+1}^{N+1} \|\Lambda^{k\alpha}u(t)\|_{L^2}^2.$$

These inequalities close the nonlinear *a priori* estimates for every $\beta \in [0, 1]$. Together with the negative Sobolev bounds, we obtain the higher-order decay estimates that also control the nonlinear remainder in the proof of sharpness of the decay exponents.

Our decay analysis follows the negative Sobolev framework introduced by Guo and Wang [13] and subsequently extended to a negative Besov framework by Sohinger and Strain [30]. For the Boltzmann equation, Guo and Wang propagated the $L_v^2\dot{H}_x^{-s}$ norms for $s \in [0, 3/2)$ and combined this control with energy estimates and interpolation to derive algebraic decay estimates. For the non-cutoff Boltzmann equation with hard and soft potentials, Sohinger and Strain established higher-order decay estimates under additional negative Besov regularity in the spatial variable. We combine the preceding energy inequalities with the propagation of negative Sobolev norms. For $0 \leq \vartheta \leq 3/2 - \alpha$, the nonlinear terms in the negative Sobolev estimates are controlled directly by the total dissipation $\mathcal{D}_N(t)$. For the remaining range $3/2 - \alpha < \vartheta < 3/2$, we first use the decay already obtained at the endpoint $\vartheta = 3/2 - \alpha$ and then prove the time integrability of the corresponding nonlinear source. Once the negative Sobolev norms are propagated for all $0 \leq \vartheta < 3/2$, interpolation between the negative and positive Sobolev levels gives (1.7). The improved density decay (1.8) follows further from the Riesz relation $\Lambda^{2\beta}\phi = \rho - 1$ and interpolation between the available potential norms. We note that the initial negative Sobolev norms in (1.5) are only assumed to be finite but need not be small.

To establish the optimality asserted in Theorem 1.1(iii), we construct low-frequency initial data along a spectral mode of the linearized system and pair conjugate profiles on opposite angular caps to obtain real-valued data. Up to smooth radial and angular cutoffs, the Fourier amplitudes of $u_0$ and $\Lambda^\beta\phi_0$ scale like $\delta|\xi|^{\vartheta+\varepsilon-3/2}$, with an additional factor $|\xi|^\beta$ for the density perturbation. The parameter $0 < \varepsilon \leq 3/2 - \vartheta$ ensures the negative Sobolev condition (1.5), while allowing the profiles to approach the integrability threshold as $\varepsilon$ decreases to zero. The low-frequency cutoff gives arbitrary positive Sobolev regularity, and, for each fixed $\varepsilon$, a sufficiently small amplitude $\delta > 0$ ensures (1.3). The associated linear evolution admits matching upper and lower decay bounds proportional to $\delta$. At the derivative orders specified in Theorem 1.1(iii), its decay exponents exceed the corresponding exponents by $\varepsilon/(2\alpha)$. We then show that the low-frequency Duhamel remainder is of order $O(\delta^2)$ and decays at least as fast as the linear contribution. For each fixed $\varepsilon$, choosing $\delta$ sufficiently small ensures that the remainder cannot cancel the leading linear contribution, thereby yielding nonlinear lower bounds with the same decay exponent. Together with the corresponding upper bounds, these estimates show that the decay dictated by the low-frequency linear dynamics persists for these nonlinear solutions, including the Riesz-induced gain of $\beta/(2\alpha)$ in the density decay exponent. Since $\varepsilon$ can be chosen arbitrarily small, these lower bounds exclude any positive improvement of the corresponding decay exponents for fixed $N$ and $\vartheta$.

The remainder of the paper is organized as follows. Section 2 establishes the energy and interactive dissipation estimates, and Section 3 develops the negative Sobolev estimates. These estimates are combined in Section 4 to prove Theorem 1.1(i)–(ii). Section 5 constructs a class of admissible initial data, proves the two-sided decay estimates in Theorem 1.6, and completes the proof of Theorem 1.1(iii). Appendix A collects the properties of the fractional Lamé–Navier operator and the analytic inequalities used in the proofs.

**Notation.** We denote the unitary Fourier transform of $f$ by $\widehat{f}$ and write $D = -\mathrm{i}\nabla$. For a scalar symbol $m$, the Fourier multiplier $m(D)$ is defined by

$$\widehat{m(D)f}(\xi) = m(\xi)\widehat{f}(\xi).$$

In particular, for $k \in \mathbb{R}$, we write $\Lambda^k = |D|^k = |\nabla|^k = (-\Delta)^{k/2}$. For convenience, we occasionally write $\nabla^k f$ for the collection of all spatial derivatives of order $k$ when $k$ is a nonnegative integer, and for $\Lambda^k f$ otherwise. We define the homogeneous Sobolev space $\dot{H}^k(\mathbb{R}^3)$ as the space of all tempered distributions $f$ for which the seminorm $\|f\|_{\dot{H}^k}$ is finite, given by

$$\|f\|_{\dot{H}^k} := \|\Lambda^k f\|_{L^2} = \||\xi|^k \widehat{f}\|_{L^2}.$$

The notation $L^p(\mathbb{R}^3)$, $1 \le p \le +\infty$, refers to the standard Lebesgue space with norm $\|\cdot\|_{L^p}$. We also denote by $(\cdot,\cdot)$ the $L^2(\mathbb{R}^3)$ inner product, with summation over components for vector-valued functions. In the estimates established in this paper, $c$ and $C$ denote generic positive constants that may vary from line to line. They may depend on fixed parameters, differentiation orders, and fixed auxiliary functions, but are independent of time and of the particular initial data within the class under consideration. Additional dependence is indicated explicitly.

## 2. Energy estimates

Letting $\varrho = \rho - 1$, we rewrite the system (1.1) in perturbation form as

$$\begin{cases} \partial_t \varrho + \mathrm{div} u = -\mathrm{div}(\varrho u), \\ \partial_t u + u\cdot\nabla u + \nabla\varrho + \nabla\phi + \mathfrak{L}^\alpha u = H(\varrho)\mathfrak{L}^\alpha u - F(\varrho)\nabla\varrho, \\ \Lambda^{2\beta}\phi = \varrho, \\ (\varrho, u)|_{t=0} = (\varrho_0, u_0), \end{cases} \tag{2.1}$$

where the two nonlinear functions of $\varrho$ are defined by

$$H(\varrho) := \frac{\varrho}{\varrho+1} \quad \text{and} \quad F(\varrho) := \frac{p'(\varrho+1)}{\varrho+1} - 1 = (\varrho+1)^{\gamma-2} - 1.$$

For subsequent estimates, we define, for $s \in \mathbb{R}$ and $u \in \dot{H}^{s+\alpha}(\mathbb{R}^3)$,

$$\mathfrak{D}_{s,L}(u) := (\mathbb{B}_\alpha \Lambda^{s+\alpha} u, \Lambda^{s+\alpha} u).$$

By (A.1), $\mathfrak{D}_{s,L}(u)$ is equivalent to $\|\Lambda^{s+\alpha}u\|_{L^2}^2$, and in particular,

$$\|\Lambda^{s+\alpha}u\|_{L^2}^2 \le \nu_-^{-1}\mathfrak{D}_{s,L}(u).$$

Hence, the Lamé dissipation $\mathfrak{D}_{s,L}(u)$ provides the coercive control of $\|\Lambda^{s+\alpha}u\|_{L^2}^2$ needed in the subsequent energy estimates. To symmetrize the hyperbolic part and thereby avoid derivative loss, we introduce the sound-speed variable

$$\mathfrak{c} := \frac{2}{\gamma-1}\left((\varrho+1)^{\frac{\gamma-1}{2}} - 1\right).$$

With this change of variable, the system (2.1) can be written equivalently as

$$\begin{cases} \partial_t \mathfrak{c} + u \cdot \nabla \mathfrak{c} + \Big(1 + \dfrac{\gamma-1}{2}\mathfrak{c}\Big)\mathrm{div} u = 0, \\ \partial_t u + u \cdot \nabla u + \Big(1 + \dfrac{\gamma-1}{2}\mathfrak{c}\Big)\nabla \mathfrak{c} + \nabla \phi + \mathfrak{L}^\alpha u = H(\varrho)\mathfrak{L}^\alpha u, \\ (-\Delta)^\beta \phi = \varrho. \end{cases} \tag{2.2}$$

Throughout this section, we assume *a priori* that, for some sufficiently small constant $\delta_1 > 0$,

$$\sqrt{\mathcal{E}_N(t)} := \sum_{k=0}^{N} \Big( \|\Lambda^{k\alpha}\varrho\|_{L^2} + \|\Lambda^{k\alpha}u\|_{L^2} + \|\Lambda^{k\alpha+\beta}\phi\|_{L^2} \Big) \le \delta_1, \tag{2.3}$$

where $N \ge 3$ is an integer satisfying $N\alpha > 5/2$. Furthermore, we define the corresponding total dissipation rate functional

$$\mathcal{D}_N(t) := \sum_{k=1}^{N} \Big( \|\Lambda^{k\alpha}\varrho\|_{L^2}^2 + \|\Lambda^{k\alpha+\beta}\phi\|_{L^2}^2 \Big) + \sum_{k=1}^{N+1} \|\Lambda^{k\alpha}u\|_{L^2}^2. \tag{2.4}$$

First, an application of the Sobolev embedding alongside (2.3) ensures that

$$\frac{3}{4} \le \rho = \varrho + 1 \le \frac{5}{4}. \tag{2.5}$$

Hence, we immediately have from Lemma A.6

$$\|H(\varrho)\|_{L^\infty} \le C\|\varrho\|_{L^\infty}, \quad \|F(\varrho)\|_{L^\infty} \le C\|\varrho\|_{L^\infty}, \tag{2.6}$$

and

$$\|\Lambda^s H(\varrho)\|_{L^p} \le C\|\Lambda^s \varrho\|_{L^p}, \quad \|\Lambda^s F(\varrho)\|_{L^p} \le C\|\Lambda^s \varrho\|_{L^p}, \tag{2.7}$$

for any $s \ge 0$ and $1 < p < +\infty$.

To transfer the ensuing energy estimates from $\mathfrak{c}$ to $\varrho$, note that (2.5) makes the transformation $\varrho \mapsto \mathfrak{c}$ a smooth diffeomorphism on the relevant ranges. Since this transformation and its inverse vanish at zero, Lemma A.6, applied in both directions, yields, for $k = 0, 1, \dots, N$,

$$c\|\Lambda^{k\alpha}\varrho\|_{L^2} \le \|\Lambda^{k\alpha}\mathfrak{c}\|_{L^2} \le C\|\Lambda^{k\alpha}\varrho\|_{L^2}. \tag{2.8}$$

With the equivalence (2.8) available, we now establish the following lower-order energy estimate.

**Lemma 2.1.** *Suppose that the* a priori *assumption* (2.3) *holds. Then, for* $k = 0, 1, \cdots, N-1$, *we have*

$$\begin{aligned} &\frac{\mathrm{d}}{\mathrm{d}t}\Big( \|\Lambda^{k\alpha}\mathfrak{c}\|_{L^2}^2 + \|\Lambda^{k\alpha}u\|_{L^2}^2 + \|\Lambda^{k\alpha+\beta}\phi\|_{L^2}^2 \Big) + \mathfrak{D}_{k\alpha,L}(u) \\ &\le C\sqrt{\mathcal{E}_N(t)}\Big( \|\Lambda^{(k+1)\alpha}\mathfrak{c}\|_{L^2}^2 + \|\Lambda^{(k+1)\alpha}u\|_{L^2}^2 + \|\Lambda^{(k+1)\alpha+\beta}\phi\|_{L^2}^2 \Big). \end{aligned} \tag{2.9}$$

*Proof.* Applying $\Lambda^{k\alpha}$ to the first two equations in (2.2), taking the $L^2$ inner product with $\Lambda^{k\alpha}\mathfrak{c}$ and $\Lambda^{k\alpha}u$, respectively, and using integration by parts, we obtain

$$\frac{1}{2}\frac{\mathrm{d}}{\mathrm{d}t}\Big( \|\Lambda^{k\alpha}\mathfrak{c}\|_{L^2}^2 + \|\Lambda^{k\alpha}u\|_{L^2}^2 \Big) + \mathfrak{D}_{k\alpha,L}(u) + (\Lambda^{k\alpha}\nabla\phi, \Lambda^{k\alpha}u) = \sum_{i=1}^{5} I_i, \tag{2.10}$$

where

$$I_1 := -\,(\Lambda^{k\alpha}(u \cdot \nabla \mathfrak{c}), \Lambda^{k\alpha}\mathfrak{c}), \qquad I_2 := -(\Lambda^{k\alpha}(u \cdot \nabla u), \Lambda^{k\alpha}u),$$

$$I_3 := -\,\frac{\gamma-1}{2}\Big(\Lambda^{k\alpha}(\mathfrak{c}\,\mathrm{div}u), \Lambda^{k\alpha}\mathfrak{c}\Big), \qquad I_4 := -\frac{\gamma-1}{2}\Big(\Lambda^{k\alpha}(\mathfrak{c}\,\nabla\mathfrak{c}), \Lambda^{k\alpha}u\Big),$$

$$I_5 := \Big(\Lambda^{k\alpha}\left(H(\varrho)\mathfrak{L}^\alpha u\right), \Lambda^{k\alpha}u\Big).$$

For $I_1$, integration by parts, Lemma A.3 and Hölder's inequality give that

$$
\begin{aligned}
I_1 &= -\left(\left[\Lambda^{k\alpha}, u\cdot\nabla\right]\mathfrak{c}, \Lambda^{k\alpha}\mathfrak{c}\right) - \left(u\cdot\nabla\Lambda^{k\alpha}\mathfrak{c}, \Lambda^{k\alpha}\mathfrak{c}\right) \\
&\le C\|\Lambda^{k\alpha}\mathfrak{c}\|_{L^{\frac{6}{3-2\alpha}}}\left(\|\Lambda^{k\alpha}u\|_{L^3}\|\nabla\mathfrak{c}\|_{L^{\frac{6}{1+2\alpha}}} + \|\nabla u\|_{L^{\frac{6}{1+2\alpha}}}\|\Lambda^{k\alpha}\mathfrak{c}\|_{L^3}\right).
\end{aligned}
\tag{2.11}
$$

With the help of the Sobolev embeddings $\dot{H}^s(\mathbb{R}^3)\hookrightarrow L^{\frac{6}{3-2s}}(\mathbb{R}^3)$ for $s\in[0,3/2)$, $\dot{H}^{\frac12}(\mathbb{R}^3)\hookrightarrow L^3(\mathbb{R}^3)$ and $\dot{H}^{1-s}(\mathbb{R}^3)\hookrightarrow L^{\frac{6}{1+2s}}(\mathbb{R}^3)$ for $s\in[0,1]$, and Lemma A.2, we have

$$
\|\Lambda^{k\alpha}\mathfrak{c}\|_{L^{\frac{6}{3-2\alpha}}} \le C\|\Lambda^{(k+1)\alpha}\mathfrak{c}\|_{L^2}, \tag{2.12}
$$

$$
\|\Lambda^{k\alpha}(u,\mathfrak{c})\|_{L^3} \le C\|\Lambda^{k\alpha+\frac12}(u,\mathfrak{c})\|_{L^2} \le C\|\Lambda^{(k+1)\alpha}(u,\mathfrak{c})\|_{L^2}^{\frac{k\alpha+\frac12}{(k+1)\alpha}}\|(u,\mathfrak{c})\|_{L^2}^{\frac{\alpha-\frac12}{(k+1)\alpha}}, \tag{2.13}
$$

$$
\|\nabla(u,\mathfrak{c})\|_{L^{\frac{6}{1+2\alpha}}} \le C\|\Lambda^{2-\alpha}(u,\mathfrak{c})\|_{L^2} \le C\|\Lambda^{(k+1)\alpha}(u,\mathfrak{c})\|_{L^2}^{\frac{\alpha-\frac12}{(k+1)\alpha}}\|\Lambda^{\theta}(u,\mathfrak{c})\|_{L^2}^{\frac{k\alpha+\frac12}{(k+1)\alpha}}, \tag{2.14}
$$

where

$$
\theta = \frac{(5-4\alpha)(k+1)\alpha}{2k\alpha+1} \in \left[\frac{N\alpha(5-4\alpha)}{2(N-1)\alpha+1}, \alpha(5-4\alpha)\right].
$$

Combining (2.12)–(2.14) with (2.11), we have from $N>5/(2\alpha)>5-4\alpha$ for $1/2\le\alpha\le1$ that

$$
I_1 \le C\sqrt{\mathcal{E}_N(t)}\left(\|\Lambda^{(k+1)\alpha}\mathfrak{c}\|_{L^2}^2 + \|\Lambda^{(k+1)\alpha}u\|_{L^2}^2\right). \tag{2.15}
$$

Applying the preceding argument with $\mathfrak{c}$ replaced by $u$, we obtain

$$
I_2 \le C\sqrt{\mathcal{E}_N(t)}\|\Lambda^{(k+1)\alpha}u\|_{L^2}^2. \tag{2.16}
$$

For $I_3$ and $I_4$, using the same argument as for $I_1$, together with integration by parts, we obtain

$$
I_3 + I_4 \le C\sqrt{\mathcal{E}_N(t)}\left(\|\Lambda^{(k+1)\alpha}\mathfrak{c}\|_{L^2}^2 + \|\Lambda^{(k+1)\alpha}u\|_{L^2}^2\right). \tag{2.17}
$$

For $I_5$, recalling from (1.2) that $\mathfrak{L}^\alpha = \mathbb{B}_\alpha\Lambda^{2\alpha}$, and using the self-adjointness of $\Lambda^\alpha$ alongside the fact that both $\mathbb{B}_\alpha$ and $\mathfrak{L}^\alpha$ commute with scalar powers of $\Lambda$, we obtain

$$
\begin{aligned}
I_5 =& \left(\left[\Lambda^{k\alpha}, H(\varrho)\right]\mathfrak{L}^\alpha u, \Lambda^{k\alpha}u\right) + \left(H(\varrho)\Lambda^{k\alpha}\mathfrak{L}^\alpha u, \Lambda^{k\alpha}u\right) \\
=& \left(\left[\Lambda^{k\alpha}, H(\varrho)\right]\mathfrak{L}^\alpha u, \Lambda^{k\alpha}u\right) + \left(\mathbb{B}_\alpha\Lambda^{(k+1)\alpha}u, \left[\Lambda^\alpha, H(\varrho)\right]\Lambda^{k\alpha}u\right) + \left(\mathbb{B}_\alpha\Lambda^{(k+1)\alpha}u, H(\varrho)\Lambda^{(k+1)\alpha}u\right) \\
:=& I_{51} + I_{52} + I_{53}.
\end{aligned}
\tag{2.18}
$$

For $I_{51}$, we only need to consider the case $1\le k\le N-1$, since this term vanishes when $k=0$. By Lemma A.3 and the $L^p$-boundedness of $\mathbb{B}_\alpha$, namely,

$$
\|\mathbb{B}_\alpha f\|_{L^p} \le C\|f\|_{L^p} \quad \text{for } 1<p<+\infty, \tag{2.19}
$$

we deduce that

$$
I_{51} \le C\|\Lambda^{k\alpha}u\|_{L^{\frac{6}{3-2\alpha}}}\left(\|\Lambda^{k\alpha}H(\varrho)\|_{L^3}\|\Lambda^{2\alpha}u\|_{L^{\frac{6}{1+2\alpha}}} + \|\nabla H(\varrho)\|_{L^3}\|\Lambda^{(k+2)\alpha-1}u\|_{L^{\frac{6}{1+2\alpha}}}\right). \tag{2.20}
$$

As in (2.13), we use (2.7) to obtain

$$
\|\Lambda^{k\alpha}H(\varrho)\|_{L^3} \le C\|\Lambda^{(k+1)\alpha}\varrho\|_{L^2}^{\frac{k\alpha+\frac12}{(k+1)\alpha}}\|\varrho\|_{L^2}^{\frac{\alpha-\frac12}{(k+1)\alpha}}.
$$

By the Sobolev embedding $\dot{H}^{1-s}(\mathbb{R}^3) \hookrightarrow L^{\frac{6}{1+2s}}(\mathbb{R}^3)$ with $s \in [0,1]$, and Lemma A.2, we have

$$\|\Lambda^{(k+2)\alpha-1}u\|_{L^{\frac{6}{1+2\alpha}}} \le C\|\Lambda^{(k+1)\alpha}u\|_{L^2},$$

$$\|\Lambda^{2\alpha}u\|_{L^{\frac{6}{1+2\alpha}}} \le C\|\Lambda^{1+\alpha}u\|_{L^2} \le C\|\Lambda^{(k+1)\alpha}u\|_{L^2}^{\frac{\alpha-\frac12}{(k+1)\alpha}}\|\Lambda^{\theta}u\|_{L^2}^{\frac{k\alpha+\frac12}{(k+1)\alpha}},$$

where

$$\theta = \frac{3(k+1)\alpha}{2k\alpha+1} \in \left[\frac{3N\alpha}{2(N-1)\alpha+1}, \frac{6\alpha}{2\alpha+1}\right].$$

Therefore, substituting the above three estimates into (2.20), and using $N \ge 3$ together with the Sobolev embeddings $\dot{H}^s(\mathbb{R}^3) \hookrightarrow L^{\frac{6}{3-2s}}(\mathbb{R}^3)$ for $s \in [0,3/2)$ and $\dot{H}^{\frac12}(\mathbb{R}^3) \hookrightarrow L^3(\mathbb{R}^3)$, we obtain

$$I_{51} \le C\sqrt{\mathcal{E}_N(t)}\left(\|\Lambda^{(k+1)\alpha}\varrho\|_{L^2}^2 + \|\Lambda^{(k+1)\alpha}u\|_{L^2}^2\right). \tag{2.21}$$

For $I_{52}$, Lemma A.3, (2.7), (2.19), the Sobolev embeddings $\dot{H}^s(\mathbb{R}^3) \hookrightarrow L^{\frac{6}{3-2s}}(\mathbb{R}^3)$ for $s \in [0,3/2)$, $\dot{H}^{\frac32-s}(\mathbb{R}^3) \hookrightarrow L^{\frac3s}(\mathbb{R}^3)$ for $s \in (0,3/2]$ and $\dot{H}^{\frac12}(\mathbb{R}^3) \hookrightarrow L^3(\mathbb{R}^3)$, together with $N\alpha > 5/2$, yield

$$\begin{aligned} I_{52} &\le C\|\Lambda^{(k+1)\alpha}u\|_{L^2}\left(\|\Lambda^{\alpha}H(\varrho)\|_{L^{\frac{3}{\alpha}}}\|\Lambda^{k\alpha}u\|_{L^{\frac{6}{3-2\alpha}}} + \|\nabla H(\varrho)\|_{L^3}\|\Lambda^{(k+1)\alpha-1}u\|_{L^6}\right) \\ &\le C\sqrt{\mathcal{E}_N(t)}\|\Lambda^{(k+1)\alpha}u\|_{L^2}^2. \end{aligned} \tag{2.22}$$

For $I_{53}$, Hölder's inequality and (2.19) imply

$$I_{53} \le C\|H(\varrho)\|_{L^\infty}\|\Lambda^{(k+1)\alpha}u\|_{L^2}^2 \le C\sqrt{\mathcal{E}_N(t)}\|\Lambda^{(k+1)\alpha}u\|_{L^2}^2. \tag{2.23}$$

Substituting (2.21), (2.22) and (2.23) into (2.18), we have

$$I_5 \le C\sqrt{\mathcal{E}_N(t)}\left(\|\Lambda^{(k+1)\alpha}\varrho\|_{L^2}^2 + \|\Lambda^{(k+1)\alpha}u\|_{L^2}^2\right). \tag{2.24}$$

It remains to handle the last term on the left-hand side of (2.10). By integration by parts, together with the first and third equations in (2.1), we obtain

$$\begin{aligned} (\Lambda^{k\alpha}\nabla\phi, \Lambda^{k\alpha}u) &= \left(\Lambda^{k\alpha}\phi, \Lambda^{k\alpha}\left(\partial_t\varrho + \operatorname{div}(\varrho u)\right)\right) \\ &= \frac12\frac{\mathrm{d}}{\mathrm{d}t}\|\Lambda^{k\alpha+\beta}\phi\|_{L^2}^2 + \left(\Lambda^{k\alpha}\phi, \Lambda^{k\alpha}\operatorname{div}(\varrho u)\right). \end{aligned} \tag{2.25}$$

Denote the last term in (2.25) by $I_6$. For $k=0$, Plancherel's theorem and the Cauchy–Schwarz inequality yield

$$|I_6| \le \|\Lambda^{\alpha+\beta}\phi\|_{L^2}\|\Lambda^{1-\alpha-\beta}(\varrho u)\|_{L^2}. \tag{2.26}$$

We first treat the endpoint $(\alpha,\beta) = (1/2, 0)$. By Lemma A.3, we have

$$\|\Lambda^{\frac12}(\varrho u)\|_{L^2} \le C\left(\|\Lambda^{\frac12}\varrho\|_{L^6}\|u\|_{L^3} + \|\Lambda^{\frac12}u\|_{L^6}\|\varrho\|_{L^3}\right) \le C\sqrt{\mathcal{E}_N(t)}(\|\Lambda^{\alpha}\varrho\|_{L^2} + \|\Lambda^{\alpha}u\|_{L^2}).$$

Away from the endpoint $(\alpha,\beta) = (1/2,0)$, we set $s_1 = 5/2 - 2\alpha - \beta < 3/2$, $s_2 = \alpha \in [1/2,1]$. Then $s_1 + s_2 = 5/2 - \alpha - \beta \in [1/2, 2]$ and $s_1 + s_2 - 3/2 = 1 - \alpha - \beta$. Hence, Lemma A.7 gives

$$\|\Lambda^{1-\alpha-\beta}(\varrho u)\|_{L^2} \le C\|\Lambda^{\frac52-2\alpha+\beta}\phi\|_{L^2}\|\Lambda^{\alpha}u\|_{L^2} \le C\sqrt{\mathcal{E}_N(t)}\|\Lambda^{\alpha}u\|_{L^2},$$

where we used $5/2 - 2\alpha < N\alpha$. Consequently, combining the above two estimates with (2.26), we have

$$|I_6| \le C\sqrt{\mathcal{E}_N(t)}\left(\|\Lambda^{\alpha}\varrho\|_{L^2}^2 + \|\Lambda^{\alpha}u\|_{L^2}^2 + \|\Lambda^{\alpha+\beta}\phi\|_{L^2}^2\right) \quad \text{for } k = 0.$$

We next consider $1 \leq k \leq N-1$. Recalling (2.25), we obtain

$$|I_6| \leq \|\Lambda^{(k+1)\alpha+\beta}\phi\|_{L^2}\|\Lambda^{1+(k-1)\alpha-\beta}(\varrho u)\|_{L^2}. \tag{2.27}$$

The choice of Sobolev exponents for estimating the second factor in (2.27) depends on the range of $\alpha+\beta$, leading to the following two cases.

**Case 1:** $\alpha+\beta \in [1/2, 3/2)$**.** In this case, we clearly have $2\alpha+\beta-1 \in [0,3/2)$. Applying Lemma A.3 and utilizing the conditions $2\alpha+\beta-1 \in [0,3/2)$ and $N\alpha > 2$, we find that

$$\begin{aligned}
&\|\Lambda^{1+(k-1)\alpha-\beta}(\varrho u)\|_{L^2} \\
&\leq C\left(\|\Lambda^{1+(k-1)\alpha-\beta}\varrho\|_{L^{\frac{6}{3-2(2\alpha+\beta-1)}}}\|u\|_{L^{\frac{3}{2\alpha+\beta-1}}} + \|\Lambda^{1+(k-1)\alpha-\beta}u\|_{L^{\frac{6}{3-2(2\alpha+\beta-1)}}}\|\varrho\|_{L^{\frac{3}{2\alpha+\beta-1}}}\right) \\
&\leq C\sqrt{\mathcal{E}_N(t)}\left(\|\Lambda^{(k+1)\alpha}\varrho\|_{L^2} + \|\Lambda^{(k+1)\alpha}u\|_{L^2}\right),
\end{aligned} \tag{2.28}$$

where we have used the Sobolev embeddings $\dot{H}^s(\mathbb{R}^3) \hookrightarrow L^{\frac{6}{3-2s}}(\mathbb{R}^3)$ for $s \in [0,3/2)$, $\dot{H}^{\frac{3}{2}-s}(\mathbb{R}^3) \hookrightarrow L^{\frac{3}{s}}(\mathbb{R}^3)$ for $s \in (0,3/2]$. At the endpoint $(\alpha,\beta) = (1/2,0)$, the $L^\infty$ factors in (2.28) are controlled by $H^{\frac{3}{2}+\delta}(\mathbb{R}^3) \hookrightarrow L^\infty(\mathbb{R}^3)$ for some $\delta > 0$, as permitted by $N\alpha > 5/2$. Substituting (2.28) into (2.27), we deduce that

$$|I_6| \leq C\sqrt{\mathcal{E}_N(t)}\left(\|\Lambda^{(k+1)\alpha}\varrho\|_{L^2}^2 + \|\Lambda^{(k+1)\alpha}u\|_{L^2}^2 + \|\Lambda^{(k+1)\alpha+\beta}\phi\|_{L^2}^2\right) \quad \text{for } k = 1,2,\cdots,N-1.$$

**Case 2:** $\alpha+\beta \in [3/2, 2]$**.** This immediately implies $1/2 \leq \beta \leq 1$. Applying Lemma A.3 yields

$$\|\Lambda^{1+(k-1)\alpha-\beta}(\varrho u)\|_{L^2} \leq C\left(\|\Lambda^{1+(k-1)\alpha-\beta}\varrho\|_{L^{\frac{6}{3-2(2\alpha-1)}}}\|u\|_{L^{\frac{3}{2\alpha-1}}} + \|\Lambda^{1+(k-1)\alpha-\beta}u\|_{L^6}\|\varrho\|_{L^3}\right). \tag{2.29}$$

Employing the Sobolev embeddings $\dot{H}^s(\mathbb{R}^3) \hookrightarrow L^{\frac{6}{3-2s}}(\mathbb{R}^3)$ for $s \in [0,3/2)$, $\dot{H}^1(\mathbb{R}^3) \hookrightarrow L^6(\mathbb{R}^3)$ and $\dot{H}^{\frac{1}{2}}(\mathbb{R}^3) \hookrightarrow L^3(\mathbb{R}^3)$, alongside Lemma A.2 and the third equation in (2.1), we have

$$\begin{aligned}
&\|\Lambda^{1+(k-1)\alpha-\beta}\varrho\|_{L^{\frac{6}{3-2(2\alpha-1)}}} \leq C\|\Lambda^{(k+1)\alpha-\beta}\varrho\|_{L^2} \leq C\|\Lambda^{(k+1)\alpha+\beta}\phi\|_{L^2}, \\
&\|\Lambda^{1+(k-1)\alpha-\beta}u\|_{L^6} \leq C\|\Lambda^{2+(k-1)\alpha-\beta}u\|_{L^2} \leq C\|\Lambda^{(k+1)\alpha}u\|_{L^2}^{\frac{2+(k-1)\alpha-\beta}{(k+1)\alpha}}\|u\|_{L^2}^{\frac{2\alpha+\beta-2}{(k+1)\alpha}}, \\
&\|\varrho\|_{L^3} \leq C\|\Lambda^{\frac{1}{2}+2\beta}\phi\|_{L^2} \leq C\|\Lambda^{(k+1)\alpha+\beta}\phi\|_{L^2}^{\frac{2\alpha+\beta-2}{(k+1)\alpha}}\|\Lambda^{\theta+\beta}\phi\|_{L^2}^{\frac{2+(k-1)\alpha-\beta}{(k+1)\alpha}},
\end{aligned}$$

where

$$\theta = \frac{\left(\frac{5}{2}-2\alpha\right)(k+1)\alpha}{2+(k-1)\alpha-\beta} \in \left[\frac{N\alpha(5-4\alpha)}{2(N-1)\alpha+1}, \alpha(5-4\alpha)\right].$$

Substituting these three estimates into (2.29) and recalling (2.27), together with $N \geq 3 \geq 5-4\alpha$ for $1/2 \leq \alpha \leq 1$, one has

$$|I_6| \leq C\sqrt{\mathcal{E}_N(t)}\left(\|\Lambda^{(k+1)\alpha}u\|_{L^2}^2 + \|\Lambda^{(k+1)\alpha+\beta}\phi\|_{L^2}^2\right) \quad \text{for } k = 1,2,\cdots,N-1.$$

Combining the estimate for $k=0$ with the two cases above, we obtain, for $0 \leq k \leq N-1$,

$$|I_6| \leq C\sqrt{\mathcal{E}_N(t)}\left(\|\Lambda^{(k+1)\alpha}\varrho\|_{L^2}^2 + \|\Lambda^{(k+1)\alpha}u\|_{L^2}^2 + \|\Lambda^{(k+1)\alpha+\beta}\phi\|_{L^2}^2\right). \tag{2.30}$$

Combining (2.10) and (2.25) with the above estimates (2.15), (2.16), (2.17), (2.24), and (2.30), and using the norm equivalence (2.8), we obtain (2.9). □

**Lemma 2.2.** *Suppose that the* a priori *assumption* (2.3) *holds. Then we have*

$$\begin{aligned}&\frac{\mathrm{d}}{\mathrm{d}t}\left(\|\Lambda^{N\alpha}\mathfrak{c}\|_{L^2}^2+\|\Lambda^{N\alpha}u\|_{L^2}^2+\left\|\frac{1}{\sqrt{\rho}}\Lambda^{N\alpha+\beta}\phi\right\|_{L^2}^2\right)+\mathfrak{D}_{N\alpha,L}(u)\\&\leq C\sqrt{\mathcal{E}_N(t)}\left(\|\Lambda^{N\alpha}\mathfrak{c}\|_{L^2}^2+\|\Lambda^{N\alpha}u\|_{L^2}^2+\|\Lambda^{(N+1)\alpha}u\|_{L^2}^2+\|\Lambda^{N\alpha+\beta}\phi\|_{L^2}^2\right).\end{aligned}\tag{2.31}$$

*Proof.* Applying $\Lambda^{N\alpha}$ to the first two equations in (2.2), taking the $L^2$ inner product with $\Lambda^{N\alpha}\mathfrak{c}$ and $\Lambda^{N\alpha}u$, respectively, and using integration by parts, we obtain

$$\frac{1}{2}\frac{\mathrm{d}}{\mathrm{d}t}\left(\|\Lambda^{N\alpha}\mathfrak{c}\|_{L^2}^2+\|\Lambda^{N\alpha}u\|_{L^2}^2\right)+\mathfrak{D}_{N\alpha,L}(u)+(\Lambda^{N\alpha}\nabla\phi,\Lambda^{N\alpha}u)=\sum_{i=1}^{5}J_i,\tag{2.32}$$

where

$$\begin{aligned}&J_1:=-\,(\Lambda^{N\alpha}(u\cdot\nabla\mathfrak{c}),\Lambda^{N\alpha}\mathfrak{c}), && J_2:=-(\Lambda^{N\alpha}(u\cdot\nabla u),\Lambda^{N\alpha}u),\\&J_3:=-\,\frac{\gamma-1}{2}(\Lambda^{N\alpha}(\mathfrak{c}\,\mathrm{div}u),\Lambda^{N\alpha}\mathfrak{c}), && J_4:=-\frac{\gamma-1}{2}(\Lambda^{N\alpha}(\mathfrak{c}\,\nabla\mathfrak{c}),\Lambda^{N\alpha}u),\\&J_5:=\left(\Lambda^{N\alpha}\left(H(\varrho)\mathfrak{L}^{\alpha}u\right),\Lambda^{N\alpha}u\right).\end{aligned}$$

For $J_1$, integration by parts, Lemma A.3, Hölder's inequality, the Sobolev embeddings $\dot{H}^s(\mathbb{R}^3)\hookrightarrow L^{\frac{6}{3-2s}}(\mathbb{R}^3)$ for $s\in[0,3/2)$, $\dot{H}^{\frac{3}{2}-s}(\mathbb{R}^3)\hookrightarrow L^{\frac{3}{s}}(\mathbb{R}^3)$ for $s\in(0,3/2]$ and $H^{\frac{3}{2}+\delta}(\mathbb{R}^3)\hookrightarrow L^{\infty}(\mathbb{R}^3)$ for $\delta>0$, together with $N\alpha>5/2$, give that

$$\begin{aligned}J_1&=-\left([\Lambda^{N\alpha},u\cdot\nabla]\mathfrak{c},\Lambda^{N\alpha}\mathfrak{c}\right)-\left(u\cdot\nabla\Lambda^{N\alpha}\mathfrak{c},\Lambda^{N\alpha}\mathfrak{c}\right)\\&\leq C\|\Lambda^{N\alpha}\mathfrak{c}\|_{L^2}\left(\|\Lambda^{N\alpha}u\|_{L^{\frac{6}{3-2\alpha}}}\|\nabla\mathfrak{c}\|_{L^{\frac{3}{\alpha}}}+\|\nabla u\|_{L^\infty}\|\Lambda^{N\alpha}\mathfrak{c}\|_{L^2}\right)\\&\leq C\sqrt{\mathcal{E}_N(t)}\left(\|\Lambda^{N\alpha}\mathfrak{c}\|_{L^2}^2+\|\Lambda^{(N+1)\alpha}u\|_{L^2}^2\right).\end{aligned}\tag{2.33}$$

For $J_2$, the same argument as for $I_1$ yields

$$J_2\leq C\sqrt{\mathcal{E}_N(t)}\|\Lambda^{(N+1)\alpha}u\|_{L^2}^2.\tag{2.34}$$

Similarly to (2.33), we deduce for $J_3$ and $J_4$ that

$$J_3+J_4\leq C\sqrt{\mathcal{E}_N(t)}\left(\|\Lambda^{N\alpha}\mathfrak{c}\|_{L^2}^2+\|\Lambda^{(N+1)\alpha}u\|_{L^2}^2\right).\tag{2.35}$$

Following the treatment of $I_5$, we again write $J_5$ as

$$\begin{aligned}J_5=&\left(\left[\Lambda^{N\alpha},H(\varrho)\right]\mathfrak{L}^{\alpha}u,\Lambda^{N\alpha}u\right)+\left(\mathbb{B}_\alpha\Lambda^{(N+1)\alpha}u,[\Lambda^{\alpha},H(\varrho)]\,\Lambda^{N\alpha}u\right)\\&+\left(\mathbb{B}_\alpha\Lambda^{(N+1)\alpha}u,H(\varrho)\Lambda^{(N+1)\alpha}u\right):=J_{51}+J_{52}+J_{53}.\end{aligned}\tag{2.36}$$

By Lemma A.3, (2.7), (2.19), Hölder's inequality, and the Sobolev embeddings $\dot{H}^s(\mathbb{R}^3)\hookrightarrow L^{\frac{6}{3-2s}}(\mathbb{R}^3)$ for $s\in[0,3/2)$, $\dot{H}^{\frac{3}{2}-s}(\mathbb{R}^3)\hookrightarrow L^{\frac{3}{s}}(\mathbb{R}^3)$ for $s\in(0,3/2]$, $\dot{H}^{\frac{1}{2}}(\mathbb{R}^3)\hookrightarrow L^3(\mathbb{R}^3)$ and $\dot{H}^{1-s}(\mathbb{R}^3)\hookrightarrow L^{\frac{6}{1+2s}}(\mathbb{R}^3)$ for $s\in[0,1]$, along with $N\alpha>5/2\geq\alpha+3/2$, we obtain the following estimate for $J_{51}$

$$\begin{aligned}J_{51}&\leq C\|\Lambda^{N\alpha}u\|_{L^{\frac{6}{3-2\alpha}}}\left(\|\Lambda^{N\alpha}H(\varrho)\|_{L^2}\|\Lambda^{2\alpha}u\|_{L^{\frac{3}{\alpha}}}+\|\nabla H(\varrho)\|_{L^3}\|\Lambda^{(N+2)\alpha-1}u\|_{L^{\frac{6}{1+2\alpha}}}\right)\\&\leq C\sqrt{\mathcal{E}_N(t)}\left(\|\Lambda^{N\alpha}\varrho\|_{L^2}^2+\|\Lambda^{(N+1)\alpha}u\|_{L^2}^2\right).\end{aligned}\tag{2.37}$$

Arguments analogous to those used in deriving (2.22) and (2.23) yield the following estimate for $J_{52}+J_{53}$

$$J_{52}+J_{53}\le C\sqrt{\mathcal{E}_N(t)}\|\Lambda^{(N+1)\alpha}u\|_{L^2}^2. \tag{2.38}$$

Substituting (2.37) and (2.38) into (2.36), we conclude that

$$J_5\le C\sqrt{\mathcal{E}_N(t)}\left(\|\Lambda^{N\alpha}\varrho\|_{L^2}^2+\|\Lambda^{(N+1)\alpha}u\|_{L^2}^2\right). \tag{2.39}$$

Let $J_6$ denote the last term on the left-hand side of (2.32). Integration by parts yields

$$J_6=-(\Lambda^{N\alpha+\beta}\phi,\Lambda^{N\alpha-\beta}\mathrm{div}u)=J_7+J_8, \tag{2.40}$$

where

$$J_7:=\left(\Lambda^{N\alpha+\beta}\phi,\frac{1}{\rho}\left[\Lambda^{N\alpha-\beta},\rho\right]\mathrm{div}u\right),\qquad J_8:=-\left(\Lambda^{N\alpha+\beta}\phi,\frac{1}{\rho}\Lambda^{N\alpha-\beta}(\rho\,\mathrm{div}u)\right).$$

For $J_7$, using (2.5), Lemma A.3 and Hölder's inequality, we have

$$|J_7|\le C\|\Lambda^{N\alpha+\beta}\phi\|_{L^2}\left(\|\Lambda^{N\alpha-\beta}\varrho\|_{L^2}\|\mathrm{div}u\|_{L^\infty}+\|\nabla\varrho\|_{L^6}\|\Lambda^{N\alpha-\beta}u\|_{L^3}\right).$$

In view of the Sobolev embeddings $\dot H^{\frac12}(\mathbb{R}^3)\hookrightarrow L^3(\mathbb{R}^3)$ and $\dot H^1(\mathbb{R}^3)\hookrightarrow L^6(\mathbb{R}^3)$, Lemma A.4 yields

$$\begin{aligned}&\|\Lambda^{N\alpha-\beta}u\|_{L^3}\le C\|\Lambda^{N\alpha-\beta+\frac12}u\|_{L^2}\le C\|\Lambda^{(N+1)\alpha}u\|_{L^2}^{\frac{N\alpha-\beta+\frac12}{(N+1)\alpha}}\|u\|_{L^2}^{\frac{\alpha+\beta-\frac12}{(N+1)\alpha}},\\ &\|\nabla\varrho\|_{L^6}\le C\|\nabla^2\varrho\|_{L^2}\le\|\Lambda^{N\alpha}\varrho\|_{L^2}^{\frac{\alpha+\beta-\frac12}{(N+1)\alpha}}\|\Lambda^\theta\varrho\|_{L^2}^{\frac{N\alpha-\beta+\frac12}{(N+1)\alpha}},\end{aligned}$$

where

$$\theta=\frac{N\alpha(\frac52-\alpha-\beta)+2\alpha}{N\alpha-\beta+\frac12}\in\left[\frac{\left[N(\frac32-\alpha)+2\right]\alpha}{N\alpha-\frac12},\frac{\left[N(\frac52-\alpha)+2\right]\alpha}{N\alpha+\frac12}\right].$$

Combining the above three estimates with the Sobolev embedding $H^{\frac32+\delta}(\mathbb{R}^3)\hookrightarrow L^\infty(\mathbb{R}^3)$ for $\delta>0$, the third equation in (2.1), and $N\alpha>5/2$, we conclude that

$$|J_7|\le C\sqrt{\mathcal{E}_N(t)}\left(\|\Lambda^{N\alpha}\varrho\|_{L^2}^2+\|\Lambda^{N\alpha+\beta}\phi\|_{L^2}^2+\|\Lambda^{(N+1)\alpha}u\|_{L^2}^2\right). \tag{2.41}$$

For $J_8$, using the first equation in (1.1) and the third equation in (2.1), we have

$$J_8=\left(\Lambda^{N\alpha+\beta}\phi,\frac{1}{\rho}\Lambda^{N\alpha-\beta}(\partial_t\rho+u\cdot\nabla\rho)\right)=\frac12\frac{\mathrm{d}}{\mathrm{d}t}\left\|\frac{1}{\sqrt\rho}\Lambda^{N\alpha+\beta}\phi\right\|_{L^2}^2+J_{81}+J_{82}+J_{83}, \tag{2.42}$$

where

$$\begin{aligned}J_{81}&:=-\frac12\left(\partial_t\left(\frac1\rho\right)\Lambda^{N\alpha+\beta}\phi,\Lambda^{N\alpha+\beta}\phi\right),\\ J_{82}&:=\left(\Lambda^{N\alpha+\beta}\phi,\frac1\rho\left[\Lambda^{N\alpha-\beta},u\cdot\nabla\right]\varrho\right),\\ J_{83}&:=\left(\Lambda^{N\alpha+\beta}\phi,\frac1\rho u\cdot\nabla\Lambda^{N\alpha-\beta}\varrho\right).\end{aligned}$$

A direct computation gives

$$J_{81}+J_{83}=-\left(\Lambda^{N\alpha+\beta}\phi,\frac1\rho\nabla\cdot u\Lambda^{N\alpha+\beta}\phi\right),$$

which, combined with Hölder's inequality, (2.5), the Sobolev embedding $H^{\frac{3}{2}+\delta}(\mathbb{R}^3) \hookrightarrow L^\infty(\mathbb{R}^3)$ for $\delta > 0$, and $N\alpha > 5/2$, implies that

$$|J_{81} + J_{83}| \leq C\|\nabla u\|_{L^\infty}\|\Lambda^{N\alpha+\beta}\phi\|_{L^2}^2 \leq C\sqrt{\mathcal{E}_N(t)}\|\Lambda^{N\alpha+\beta}\phi\|_{L^2}^2. \tag{2.43}$$

As for $J_{82}$, an argument similar to that for $J_7$ gives

$$|J_{82}| \leq C\sqrt{\mathcal{E}_N(t)}\left(\|\Lambda^{N\alpha}\varrho\|_{L^2}^2 + \|\Lambda^{N\alpha+\beta}\phi\|_{L^2}^2 + \|\Lambda^{(N+1)\alpha}u\|_{L^2}^2\right). \tag{2.44}$$

Combining (2.42), (2.43), and (2.44), we obtain

$$J_8 \geq \frac{1}{2}\frac{\mathrm{d}}{\mathrm{d}t}\left\|\frac{1}{\sqrt{\rho}}\Lambda^{N\alpha+\beta}\phi\right\|_{L^2}^2 - C\sqrt{\mathcal{E}_N(t)}\left(\|\Lambda^{N\alpha}\varrho\|_{L^2}^2 + \|\Lambda^{N\alpha+\beta}\phi\|_{L^2}^2 + \|\Lambda^{(N+1)\alpha}u\|_{L^2}^2\right). \tag{2.45}$$

Consequently, substituting (2.41) and (2.45) into (2.40), we arrive at

$$J_6 \geq \frac{1}{2}\frac{\mathrm{d}}{\mathrm{d}t}\left\|\frac{1}{\sqrt{\rho}}\Lambda^{N\alpha+\beta}\phi\right\|_{L^2}^2 - C\sqrt{\mathcal{E}_N(t)}\left(\|\Lambda^{N\alpha}\varrho\|_{L^2}^2 + \|\Lambda^{N\alpha+\beta}\phi\|_{L^2}^2 + \|\Lambda^{(N+1)\alpha}u\|_{L^2}^2\right). \tag{2.46}$$

Finally, inserting (2.33), (2.34), (2.35), (2.39), and (2.46) into (2.32), and applying (2.8), we obtain (2.31). □

The following interactive energy estimate recovers the dissipation of $\varrho$ and $\Lambda^\beta\phi$.

**Lemma 2.3.** *Suppose that the* a priori *assumption* (2.3) *holds. Then for* $k = 0, 1, \ldots, N-1$, *we have*

$$\begin{aligned}&2\frac{\mathrm{d}}{\mathrm{d}t}\int_{\mathbb{R}^3}\Lambda^{k\alpha}u\cdot\nabla\Lambda^{k\alpha+2\alpha-2}\varrho\,\mathrm{d}x + \|\Lambda^{(k+1)\alpha}\varrho\|_{L^2}^2 + \|\Lambda^{(k+1)\alpha+\beta}\phi\|_{L^2}^2\\ &\leq C\left(\|\Lambda^{(k+1)\alpha}u\|_{L^2}^2 + \|\Lambda^{(k+2)\alpha}u\|_{L^2}^2\right).\end{aligned} \tag{2.47}$$

*Proof.* Applying $\Lambda^{k\alpha}$ to the second equation in (2.1), taking the $L^2$ inner product with $\nabla\Lambda^{k\alpha+2\alpha-2}\varrho$, and using the third equation in (2.1) together with integration by parts, we obtain

$$\frac{\mathrm{d}}{\mathrm{d}t}\int_{\mathbb{R}^3}\Lambda^{k\alpha}u\cdot\nabla\Lambda^{k\alpha+2\alpha-2}\varrho\,\mathrm{d}x + \|\Lambda^{(k+1)\alpha}\varrho\|_{L^2}^2 + \|\Lambda^{(k+1)\alpha+\beta}\phi\|_{L^2}^2 = \sum_{i=1}^{5}K_i, \tag{2.48}$$

where

$$\begin{aligned}K_1 &:= \left(\Lambda^{k\alpha}u, \nabla\Lambda^{k\alpha+2\alpha-2}\partial_t\varrho\right), & K_2 &:= -\left(\Lambda^{k\alpha}\mathfrak{L}^\alpha u, \nabla\Lambda^{k\alpha+2\alpha-2}\varrho\right),\\ K_3 &:= -\left(\Lambda^{k\alpha}(u\cdot\nabla u), \nabla\Lambda^{k\alpha+2\alpha-2}\varrho\right), & K_4 &:= \left(\Lambda^{k\alpha}(H(\varrho)\mathfrak{L}^\alpha u), \nabla\Lambda^{k\alpha+2\alpha-2}\varrho\right),\\ K_5 &:= -\left(\Lambda^{k\alpha}(F(\varrho)\nabla\varrho), \nabla\Lambda^{k\alpha+2\alpha-2}\varrho\right).\end{aligned}$$

For $K_1$, we utilize the first equation in (2.1), integration by parts, Hölder's inequality, Lemma A.3, the Sobolev embedding $H^{\frac{3}{2}+\delta}(\mathbb{R}^3) \hookrightarrow L^\infty(\mathbb{R}^3)$ for $\delta > 0$ and $N\alpha > 5/2$ to yield

$$\begin{aligned}K_1 =& -\left(\Lambda^{k\alpha}u, \nabla\Lambda^{k\alpha+2\alpha-2}(\mathrm{div}u + \mathrm{div}(\varrho u))\right)\\ &\leq \|\Lambda^{(k+1)\alpha}u\|_{L^2}^2 + \|\Lambda^{(k+1)\alpha}u\|_{L^2}\|\Lambda^{(k+1)\alpha}(\varrho u)\|_{L^2}\\ &\leq \|\Lambda^{(k+1)\alpha}u\|_{L^2}^2 + C\sqrt{\mathcal{E}_N(t)}\left(\|\Lambda^{(k+1)\alpha}u\|_{L^2}^2 + \|\Lambda^{(k+1)\alpha}\varrho\|_{L^2}^2\right).\end{aligned} \tag{2.49}$$

Since $\nabla\Lambda^{k\alpha+2\alpha-2}\varrho$ is a gradient field, its $L^2$ pairing with the transverse component of $u$ vanishes. Consequently, only the longitudinal viscous term contributes to $K_2$, and

$$K_2 = -\nu_{\mathrm{L}}\left(\Lambda^{(k+2)\alpha}u, \nabla\Lambda^{k\alpha+2\alpha-2}\varrho\right).$$

With the help of integration by parts, Hölder's inequality, the fact that $1/2 \le \alpha \le 1$ implies $0 \le 2\alpha - 1 \le \alpha$ and the Cauchy–Schwarz inequality, we have

$$
\begin{aligned}
K_2 &\le C\|\Lambda^{(k+1)\alpha+2\alpha-1}u\|_{L^2}\|\Lambda^{(k+1)\alpha}\varrho\|_{L^2} \\
&\le \eta\|\Lambda^{(k+1)\alpha}\varrho\|_{L^2}^2 + C_\eta\left(\|\Lambda^{(k+1)\alpha}u\|_{L^2}^2 + \|\Lambda^{(k+2)\alpha}u\|_{L^2}^2\right).
\end{aligned} \tag{2.50}
$$

For $K_3$, applying Hölder's inequality, Lemma A.3 and the Sobolev embeddings $\dot{H}^{1-s}(\mathbb{R}^3) \hookrightarrow L^{\frac{6}{1+2s}}(\mathbb{R}^3)$ for $s \in [0,1]$, $\dot{H}^s(\mathbb{R}^3) \hookrightarrow L^{\frac{6}{3-2s}}(\mathbb{R}^3)$ for $s \in [0,3/2)$, and $\dot{H}^{\frac{3}{2}-s}(\mathbb{R}^3) \hookrightarrow L^{\frac{3}{s}}(\mathbb{R}^3)$ for $s \in (0,3/2]$, we have

$$
\begin{aligned}
K_3 &\le C\|\nabla\Lambda^{k\alpha+2\alpha-2}\varrho\|_{L^{\frac{6}{1+2\alpha}}}\left(\|\Lambda^{k\alpha}u\|_{L^{\frac{6}{3-2\alpha}}}\|\nabla u\|_{L^3} + \|u\|_{L^{\frac{3}{1-\alpha}}}\|\Lambda^{k\alpha+1}u\|_{L^2}\right) \\
&\le C\sqrt{\mathcal{E}_N(t)}\left(\|\Lambda^{(k+1)\alpha}\varrho\|_{L^2}^2 + \|\Lambda^{(k+1)\alpha}u\|_{L^2}^2 + \|\Lambda^{(k+2)\alpha}u\|_{L^2}^2\right).
\end{aligned} \tag{2.51}
$$

Here, we have used the fact that $\alpha \le 1 \le 2\alpha$ and Lemma A.4, which imply that

$$
\|\Lambda^{k\alpha+1}u\|_{L^2} \le C\left(\|\Lambda^{(k+1)\alpha}u\|_{L^2} + \|\Lambda^{(k+2)\alpha}u\|_{L^2}\right).
$$

Next, for $K_4$, the same argument as for $K_3$, together with (2.19), yields

$$
K_4 \le C\sqrt{\mathcal{E}_N(t)}\left(\|\Lambda^{(k+1)\alpha}\varrho\|_{L^2}^2 + \|\Lambda^{(k+2)\alpha}u\|_{L^2}^2\right). \tag{2.52}
$$

For $K_5$, a straightforward calculation shows that

$$
F(\varrho)\nabla\varrho = \nabla\mathbf{F}(\varrho),
$$

where

$$
\mathbf{F}(\varrho) = (\gamma-2)\varrho^2\int_0^1(1-\theta)(1+\theta\varrho)^{\gamma-3}\,\mathrm{d}\theta := \varrho\mathcal{F}(\varrho).
$$

Clearly, $\mathcal{F}(0) = 0$ and $\mathcal{F}(\varrho)$ is smooth on the range of $\varrho$. Hence, performing integration by parts and employing Lemma A.3, Lemma A.6, the Sobolev embedding $H^{\frac{3}{2}+\delta}(\mathbb{R}^3) \hookrightarrow L^\infty(\mathbb{R}^3)$ for $\delta > 0$, and $N\alpha > 5/2$, we deduce that

$$
\begin{aligned}
K_5 =& -\left(\Lambda^{(k+1)\alpha}(\varrho\mathcal{F}(\varrho)), \Lambda^{(k+1)\alpha}\varrho\right) \\
\le& \, C\|\Lambda^{(k+1)\alpha}\varrho\|_{L^2}\left(\|\Lambda^{(k+1)\alpha}\varrho\|_{L^2}\|\mathcal{F}(\varrho)\|_{L^\infty} + \|\varrho\|_{L^\infty}\|\Lambda^{(k+1)\alpha}(\mathcal{F}(\varrho))\|_{L^2}\right) \\
\le& \, C\sqrt{\mathcal{E}_N(t)}\|\Lambda^{(k+1)\alpha}\varrho\|_{L^2}^2.
\end{aligned} \tag{2.53}
$$

Finally, substituting the estimates (2.49)–(2.53) into (2.48), and choosing $\eta > 0$ sufficiently small, we utilize the *a priori* assumption (2.3) to arrive at the desired estimate (2.47). □

## 3. Negative Sobolev estimates

In this section, we establish the following differential inequality for propagating the negative Sobolev norms of solutions to (2.1).

**Lemma 3.1.** *Suppose that the* a priori *assumption* (2.3) *holds. Then, for any* $\vartheta \in (0,3/2)$, *we have*

$$
\frac{\mathrm{d}}{\mathrm{d}t}\mathcal{E}_{-\vartheta}(t) + \mathfrak{D}_{-\vartheta,L}(u) \le C\mathcal{N}(\vartheta)\sqrt{\mathcal{E}_{-\vartheta}(t)}, \tag{3.1}
$$

*where the negative-order Sobolev energy functional is defined as*

$$
\mathcal{E}_{-\vartheta}(t) := \|\Lambda^{-\vartheta}\varrho(t)\|_{L^2}^2 + \|\Lambda^{-\vartheta}u(t)\|_{L^2}^2 + \|\Lambda^{-\vartheta+\beta}\phi(t)\|_{L^2}^2, \tag{3.2}
$$

*and*

$$\mathcal{N}(\vartheta) := \begin{cases} \mathcal{D}_N(t), & \vartheta \in (0, \frac{3}{2}-\alpha], \\ \|(\varrho, u, \Lambda^\beta \phi)\|_{L^2}^{1-\frac{3-2\vartheta}{2\alpha}} \left(\|\Lambda^\alpha(\varrho, u, \Lambda^\beta\phi)\|_{L^2} + \|\Lambda^{2\alpha}(\varrho, u)\|_{L^2}\right)^{1+\frac{3-2\vartheta}{2\alpha}}, & \vartheta \in (\frac{3}{2}-\alpha, \frac{3}{2}). \end{cases} \tag{3.3}$$

*Proof.* Applying $\Lambda^{-\vartheta}$ to the first two equations in (2.1), taking the $L^2$ inner product with $\Lambda^{-\vartheta}\varrho$ and $\Lambda^{-\vartheta}u$, respectively, and performing integration by parts, we obtain

$$\frac{1}{2}\frac{\mathrm{d}}{\mathrm{d}t}\left(\|\Lambda^{-\vartheta}\varrho\|_{L^2}^2 + \|\Lambda^{-\vartheta}u\|_{L^2}^2\right) + \mathfrak{D}_{-\vartheta,L}(u) + (\Lambda^{-\vartheta}\nabla\phi, \Lambda^{-\vartheta}u) = \sum_{i=1}^{5} L_i, \tag{3.4}$$

where

$$L_1 := -(\Lambda^{-\vartheta}(\varrho\,\mathrm{div}u), \Lambda^{-\vartheta}\varrho), \quad L_2 := -(\Lambda^{-\vartheta}(u\cdot\nabla\varrho), \Lambda^{-\vartheta}\varrho), \quad L_3 := -(\Lambda^{-\vartheta}(u\cdot\nabla u), \Lambda^{-\vartheta}u),$$
$$L_4 := (\Lambda^{-\vartheta}(H(\varrho)\mathfrak{L}^\alpha u), \Lambda^{-\vartheta}u), \quad L_5 := -(\Lambda^{-\vartheta}(F(\varrho)\nabla\varrho), \Lambda^{-\vartheta}u).$$

For $L_1$, applying Hölder's inequality and Lemma A.7 with $s_1 = 3/2-\vartheta$ and $s_2 = 0$, we obtain

$$L_1 \le \|\Lambda^{-\vartheta}(\varrho\,\mathrm{div}u)\|_{L^2}\|\Lambda^{-\vartheta}\varrho\|_{L^2} \le C\|\Lambda^{\frac{3}{2}-\vartheta}\varrho\|_{L^2}\|\nabla u\|_{L^2}\|\Lambda^{-\vartheta}\varrho\|_{L^2}. \tag{3.5}$$

The assumptions of Lemma A.7 are satisfied here since $0 < \vartheta < 3/2$, and hence $3/2-\vartheta \in (0, 3/2)$. Similarly, for $L_2$–$L_5$, we have

$$\begin{aligned} \sum_{i=2,3,5} L_i &\le C\left(\|\Lambda^{\frac{3}{2}-\vartheta}\varrho\|_{L^2} + \|\Lambda^{\frac{3}{2}-\vartheta}u\|_{L^2}\right)(\|\nabla\varrho\|_{L^2} + \|\nabla u\|_{L^2})\left(\|\Lambda^{-\vartheta}\varrho\|_{L^2} + \|\Lambda^{-\vartheta}u\|_{L^2}\right), \\ L_4 &\le C\|\Lambda^{\frac{3}{2}-\vartheta}\varrho\|_{L^2}\|\Lambda^{2\alpha}u\|_{L^2}\|\Lambda^{-\vartheta}u\|_{L^2}. \end{aligned} \tag{3.6}$$

Since the control of $\|\Lambda^{\frac{3}{2}-\vartheta}(\varrho, u)\|_{L^2}$ depends on whether $\frac{3}{2}-\vartheta \ge \alpha$, we distinguish the following two cases.

**Case 1:** $\vartheta \in (0, 3/2-\alpha]$**.** In this case, we have $\alpha \le \frac{3}{2}-\vartheta \le N\alpha$. Hence, by the definition of $\mathcal{D}_N(t)$ in (2.4), Lemma A.4, and the fact that $\alpha \le 1 \le 2\alpha$, we deduce

$$\|\Lambda^{\frac{3}{2}-\vartheta}\varrho\|_{L^2} + \|\Lambda^{\frac{3}{2}-\vartheta}u\|_{L^2} + \|\nabla\varrho\|_{L^2} + \|\nabla u\|_{L^2} + \|\Lambda^{2\alpha}u\|_{L^2} \le C\sqrt{\mathcal{D}_N(t)}. \tag{3.7}$$

Therefore, from (3.5)–(3.6), we get

$$\sum_{i=1}^{5} L_i \le C\mathcal{D}_N(t)\left(\|\Lambda^{-\vartheta}\varrho\|_{L^2} + \|\Lambda^{-\vartheta}u\|_{L^2}\right). \tag{3.8}$$

**Case 2:** $\vartheta \in (3/2-\alpha, 3/2)$**.** In this case, we have $0 < 3/2-\vartheta < \alpha$, and the derivative order $3/2-\vartheta$ cannot be directly controlled by the dissipation. Thus we use Lemma A.4 to estimate

$$\begin{aligned} &\|\Lambda^{\frac{3}{2}-\vartheta}f\|_{L^2} \le C\|f\|_{L^2}^{1-\frac{3-2\vartheta}{2\alpha}}\|\Lambda^\alpha f\|_{L^2}^{\frac{3-2\vartheta}{2\alpha}}, \\ &\|\nabla f\|_{L^2} \le C\|\Lambda^\alpha f\|_{L^2}^{2-\frac{1}{\alpha}}\|\Lambda^{2\alpha}f\|_{L^2}^{\frac{1}{\alpha}-1} \le C\left(\|\Lambda^\alpha f\|_{L^2} + \|\Lambda^{2\alpha}f\|_{L^2}\right). \end{aligned} \tag{3.9}$$

Consequently, we have

$$\sum_{i=1}^{5} L_i \le C\|(\varrho, u)\|_{L^2}^{1-\frac{3-2\vartheta}{2\alpha}}\left(\|\Lambda^\alpha(\varrho, u)\|_{L^2} + \|\Lambda^{2\alpha}(\varrho, u)\|_{L^2}\right)^{1+\frac{3-2\vartheta}{2\alpha}}\left(\|\Lambda^{-\vartheta}\varrho\|_{L^2} + \|\Lambda^{-\vartheta}u\|_{L^2}\right). \tag{3.10}$$

Finally, we address the last term on the left-hand side of (3.4). Using the first and third equations in (2.1), together with integration by parts, one has

$$\begin{aligned}(\Lambda^{-\vartheta}\nabla\phi,\Lambda^{-\vartheta}u) &= -(\Lambda^{-\vartheta}\phi,\Lambda^{-\vartheta}\mathrm{div}u)\\ &= (\Lambda^{-\vartheta}\phi,\Lambda^{-\vartheta}\partial_t\varrho)+\Big(\Lambda^{-\vartheta}\phi,\Lambda^{-\vartheta}\mathrm{div}(\varrho u)\Big)\\ &= \frac{1}{2}\frac{\mathrm{d}}{\mathrm{d}t}\|\Lambda^{-\vartheta+\beta}\phi\|_{L^2}^2+\Big(\Lambda^{-\vartheta+\beta}\phi,\Lambda^{-\vartheta-\beta}\mathrm{div}(\varrho u)\Big).\end{aligned}\tag{3.11}$$

It remains to estimate the last term on the right-hand side of (3.11). Applying Hölder's inequality and Lemma A.7 with $s_1 = 3/2-\vartheta-\beta$ and $s_2 = 1$, we deduce that for any $\vartheta\in(0,3/2)$,

$$\begin{aligned}\Big|\Big(\Lambda^{-\vartheta+\beta}\phi,\Lambda^{-\vartheta-\beta}\mathrm{div}(\varrho u)\Big)\Big| &\le C\|\Lambda^{-\vartheta+\beta}\phi\|_{L^2}\|\Lambda^{\frac{3}{2}-\vartheta-\beta}\varrho\|_{L^2}\|\Lambda u\|_{L^2}\\ &\le C\|\Lambda^{-\vartheta+\beta}\phi\|_{L^2}\|\Lambda^{\frac{3}{2}-\vartheta+\beta}\phi\|_{L^2}\|\Lambda u\|_{L^2}.\end{aligned}\tag{3.12}$$

We estimate the right-hand side of (3.12) separately in the following two cases.
**Case 1:** $\vartheta\in(0,3/2-\alpha]$**.** Using the same argument as in (3.7), we immediately obtain

$$\Big|\Big(\Lambda^{-\vartheta+\beta}\phi,\Lambda^{-\vartheta-\beta}\mathrm{div}(\varrho u)\Big)\Big|\le C\mathcal{D}_N(t)\|\Lambda^{-\vartheta+\beta}\phi\|_{L^2}.\tag{3.13}$$

**Case 2:** $\vartheta\in(3/2-\alpha,3/2)$**.** Combining (3.12) with (3.9), we have

$$\begin{aligned}&\Big|\Big(\Lambda^{-\vartheta+\beta}\phi,\Lambda^{-\vartheta-\beta}\mathrm{div}(\varrho u)\Big)\Big|\\ &\le C\|\Lambda^{\beta}\phi\|_{L^2}^{1-\frac{3-2\vartheta}{2\alpha}}\Big(\|\Lambda^{\alpha}(\Lambda^{\beta}\phi,u)\|_{L^2}+\|\Lambda^{2\alpha}u\|_{L^2}\Big)^{1+\frac{3-2\vartheta}{2\alpha}}\|\Lambda^{-\vartheta+\beta}\phi\|_{L^2}.\end{aligned}\tag{3.14}$$

Finally, combining (3.4), (3.8), (3.10), (3.11), (3.13), and (3.14), we obtain (3.1). □

## 4. Global well-posedness and decay estimates

In this section, we prove Theorem 1.1(i)–(ii). The optimality assertion in Theorem 1.1(iii) is proved in Section 5.5. We first invoke the standard local well-posedness theory for smooth solutions, following the framework of Wang and Zhang [37]. For the additional parameter regimes considered here, the estimates established in Section 2 provide the required control of the electric coupling. Thus, it suffices to establish the uniform *a priori* bound required by the continuation argument.

*Proof of global well-posedness and* (1.4). Let $N\ge 3$ satisfy $N\alpha>5/2$, and let $l=0,\dots,N-1$. We close the estimates at the $l\alpha$-th spatial derivative level by defining

$$\begin{aligned}\mathcal{E}_N^l(t) := &\sum_{k=l}^{N}\Big(\|\Lambda^{k\alpha}\mathfrak{c}\|_{L^2}^2+\|\Lambda^{k\alpha}u\|_{L^2}^2\Big)+\sum_{k=l}^{N-1}\|\Lambda^{k\alpha+\beta}\phi\|_{L^2}^2+\left\|\frac{1}{\sqrt{\rho}}\Lambda^{N\alpha+\beta}\phi\right\|_{L^2}^2\\ &+2\kappa\sum_{k=l}^{N-1}\int_{\mathbb{R}^3}\Lambda^{k\alpha}u\cdot\nabla\Lambda^{k\alpha+2\alpha-2}\varrho\,\mathrm{d}x,\end{aligned}$$

where $0<\kappa\ll 1$ will be chosen sufficiently small. By (2.8), (2.5) and the Cauchy–Schwarz inequality, we have

$$c\mathcal{E}_N^l(t)\le\sum_{k=l}^{N}\Big(\|\Lambda^{k\alpha}\varrho\|_{L^2}^2+\|\Lambda^{k\alpha}u\|_{L^2}^2+\|\Lambda^{k\alpha+\beta}\phi\|_{L^2}^2\Big)\le C\mathcal{E}_N^l(t).\tag{4.1}$$

We also define the corresponding dissipation functional by

$$\mathcal{D}_N^l(t) := \sum_{k=l+1}^{N+1}\|\Lambda^{k\alpha}u\|_{L^2}^2+\sum_{k=l+1}^{N}\Big(\|\Lambda^{k\alpha}\varrho\|_{L^2}^2+\|\Lambda^{k\alpha+\beta}\phi\|_{L^2}^2\Big).$$

In light of the definitions of $\mathcal{E}_N(t)$ and $\mathcal{D}_N(t)$ in (2.3) and (2.4), respectively, there are constants $c, C > 0$ such that

$$c\mathcal{E}_N^0(t) \le \mathcal{E}_N(t) \le C\mathcal{E}_N^0(t), \qquad \mathcal{D}_N(t) = \mathcal{D}_N^0(t).$$

Summing (2.9) in Lemma 2.1 over $k = l, \dots, N-1$, adding (2.31) in Lemma 2.2, and using (A.1), (2.8), and the *a priori* assumption (2.3), we obtain

$$\begin{aligned}&\frac{\mathrm{d}}{\mathrm{d}t}\left\{\sum_{k=l}^{N}\left(\|\Lambda^{k\alpha}\mathfrak{c}\|_{L^2}^2 + \|\Lambda^{k\alpha}u\|_{L^2}^2\right) + \sum_{k=l}^{N-1}\|\Lambda^{k\alpha+\beta}\phi\|_{L^2}^2 + \left\|\frac{1}{\sqrt{\rho}}\Lambda^{N\alpha+\beta}\phi\right\|_{L^2}^2\right\}\\ &+\frac{\nu_-}{2}\sum_{k=l+1}^{N+1}\|\Lambda^{k\alpha}u\|_{L^2}^2 \le C\delta_1\sum_{k=l+1}^{N}\left(\|\Lambda^{k\alpha}\varrho\|_{L^2}^2 + \|\Lambda^{k\alpha+\beta}\phi\|_{L^2}^2\right).\end{aligned} \tag{4.2}$$

On the other hand, summing (2.47) in Lemma 2.3 over $k = l, \dots, N-1$ gives

$$2\frac{\mathrm{d}}{\mathrm{d}t}\sum_{k=l}^{N-1}\int_{\mathbb{R}^3}\Lambda^{k\alpha}u\cdot\nabla\Lambda^{k\alpha+2\alpha-2}\varrho\,\mathrm{d}x + \sum_{k=l+1}^{N}\left(\|\Lambda^{k\alpha}\varrho\|_{L^2}^2 + \|\Lambda^{k\alpha+\beta}\phi\|_{L^2}^2\right) \le C\sum_{k=l+1}^{N+1}\|\Lambda^{k\alpha}u\|_{L^2}^2. \tag{4.3}$$

We first choose $\kappa > 0$ sufficiently small such that $C\kappa \le \nu_-/4$, and then choose $\delta_1 > 0$ small enough so that $C\delta_1 \le \kappa/2$. Multiplying (4.3) by $\kappa$ and adding the result to (4.2), we obtain that for $l = 0, \dots, N-1$,

$$\frac{\mathrm{d}}{\mathrm{d}t}\mathcal{E}_N^l(t) + c\mathcal{D}_N^l(t) \le 0. \tag{4.4}$$

In particular, taking $l = 0$ and integrating directly in time gives

$$\mathcal{E}_N(t) + \int_0^t \mathcal{D}_N(\tau)\,\mathrm{d}\tau \le C\mathcal{E}_N(0),$$

which yields (1.4). Moreover, a standard continuity argument closes the *a priori* estimate (2.3), provided that $\delta_0$ in (1.3) is sufficiently small.

□

We next treat the case in which the negative Sobolev estimate can be integrated directly using the global energy bound.

*Proof of* (1.6)–(1.8) *for* $\vartheta \in [0, 3/2-\alpha]$. For the case $\vartheta = 0$, the estimate (1.6) follows directly from (1.4). We therefore only need to consider the case $\vartheta \in (0, 3/2-\alpha]$ for the negative Sobolev estimate.

Integrating (3.1) in time and utilizing the estimate (1.4), we obtain that for $\vartheta \in (0, 3/2-\alpha]$,

$$\mathcal{E}_{-\vartheta}(t) \le \mathcal{E}_{-\vartheta}(0) + C\int_0^t \mathcal{D}_N(\tau)\sqrt{\mathcal{E}_{-\vartheta}(\tau)}\,\mathrm{d}\tau \le CM_0^2 + C\delta_0^2\sup_{0\le\tau\le t}\sqrt{\mathcal{E}_{-\vartheta}(\tau)}.$$

In view of the definition (3.2), this, together with the case $\vartheta = 0$ obtained from (1.4), implies (1.6) for $\vartheta \in [0, 3/2-\alpha]$

$$\|\Lambda^{-\vartheta}\varrho(t)\|_{L^2} + \|\Lambda^{-\vartheta}u(t)\|_{L^2} + \|\Lambda^{-\vartheta+\beta}\phi(t)\|_{L^2} \le C\left(M_0 + \delta_0\right). \tag{4.5}$$

We now derive the decay estimate. For $l = 0, \dots, N-1$ with $l\alpha + \vartheta > 0$, Lemma A.4 yields

$$\|\Lambda^{(l+1)\alpha}f\|_{L^2} \ge C\|\Lambda^{-\vartheta}f\|_{L^2}^{-\frac{\alpha}{l\alpha+\vartheta}}\|\Lambda^{l\alpha}f\|_{L^2}^{1+\frac{\alpha}{l\alpha+\vartheta}}.$$

Combining this interpolation inequality with (4.5), we find

$$\begin{aligned}&\|\Lambda^{(l+1)\alpha}\varrho\|_{L^2}^2 + \|\Lambda^{(l+1)\alpha}u\|_{L^2}^2 + \|\Lambda^{(l+1)\alpha+\beta}\phi\|_{L^2}^2\\ &\ge C\left(M_0+\delta_0\right)^{-\frac{2\alpha}{l\alpha+\vartheta}}\left(\|\Lambda^{l\alpha}\varrho\|_{L^2}^2 + \|\Lambda^{l\alpha}u\|_{L^2}^2 + \|\Lambda^{l\alpha+\beta}\phi\|_{L^2}^2\right)^{1+\frac{\alpha}{l\alpha+\vartheta}}.\end{aligned}$$

This, together with (1.4), implies in particular that for $l=0,\dots,N-1$ with $l\alpha+\vartheta>0$,

$$\mathcal{D}_N^l(t) \geq C\,(M_0+\delta_0)^{-\frac{2\alpha}{l\alpha+\vartheta}}\,(\mathcal{E}_N^l(t))^{1+\frac{\alpha}{l\alpha+\vartheta}}.$$

Inserting this estimate into (4.4) yields

$$\frac{\mathrm{d}}{\mathrm{d}t}\mathcal{E}_N^l(t) + C\,(M_0+\delta_0)^{-\frac{2\alpha}{l\alpha+\vartheta}}\left(\mathcal{E}_N^l(t)\right)^{1+\frac{\alpha}{l\alpha+\vartheta}} \leq 0 \quad \text{for } \; l=0,\dots,N-1 \; \text{ with } \; l\alpha+\vartheta>0.$$

Solving this differential inequality directly gives

$$\mathcal{E}_N^l(t) \leq C\,(M_0+\delta_0)^2\,(1+t)^{-\frac{l\alpha+\vartheta}{\alpha}} \quad \text{for } \; l=0,\dots,N-1 \; \text{ with } \; l\alpha+\vartheta>0. \tag{4.6}$$

For the only excluded case $(l,\vartheta)=(0,0)$, (1.4) gives $\mathcal{E}_N(t)\leq C\delta_0^2$, which is exactly (1.7) with decay exponent zero. This, together with (4.6), establishes that for $\vartheta\in[0,3/2-\alpha]$ and $l=0,\dots,N-1$,

$$\sum_{k=l}^{N}\left(\|\Lambda^{k\alpha}\varrho(t)\|_{L^2}+\|\Lambda^{k\alpha}u(t)\|_{L^2}+\|\Lambda^{k\alpha+\beta}\phi(t)\|_{L^2}\right) \leq C\,(M_0+\delta_0)\,(1+t)^{-\frac{\vartheta}{2\alpha}-\frac{l}{2}}. \tag{4.7}$$

Next, we derive the improved decay rate for $\varrho$. For $l=0,\dots,N-2$ and $\beta\leq\alpha$, we have from the third equation in (2.1) and (4.7)

$$\|\Lambda^{l\alpha}\varrho\|_{L^2} = \|\Lambda^{l\alpha+2\beta}\phi\|_{L^2} \leq C\|\Lambda^{l\alpha+\beta}\phi\|_{L^2}^{1-\frac{\beta}{\alpha}}\|\Lambda^{(l+1)\alpha+\beta}\phi\|_{L^2}^{\frac{\beta}{\alpha}} \leq C\,(M_0+\delta_0)\,(1+t)^{-\frac{\vartheta}{2\alpha}-\frac{\beta}{2\alpha}-\frac{l}{2}}.$$

For $l=0,\dots,N-3$ and $\beta>\alpha$, we have from the third equation in (2.1), $\beta\leq 1\leq 2\alpha$ and (4.7)

$$\|\Lambda^{l\alpha}\varrho(t)\|_{L^2} \leq C\|\Lambda^{(l+1)\alpha+\beta}\phi(t)\|_{L^2}^{2-\frac{\beta}{\alpha}}\|\Lambda^{(l+2)\alpha+\beta}\phi(t)\|_{L^2}^{\frac{\beta}{\alpha}-1} \leq C\,(M_0+\delta_0)\,(1+t)^{-\frac{\vartheta}{2\alpha}-\frac{\beta}{2\alpha}-\frac{l}{2}}.$$

As for $l=N-2$ and $\beta>\alpha$, we have

$$\|\Lambda^{(N-2)\alpha}\varrho(t)\|_{L^2} \leq C\|\Lambda^{(N-1)\alpha+\beta}\phi(t)\|_{L^2}^{2-\frac{\beta}{\alpha}}\|\Lambda^{N\alpha+\beta}\phi(t)\|_{L^2}^{\frac{\beta}{\alpha}-1} \leq C\,(M_0+\delta_0)\,(1+t)^{-\frac{\vartheta}{2\alpha}-\frac{N-1}{2}}. \tag{4.8}$$

Hence, combining (4.7)–(4.8) and Sobolev interpolation, we get (1.7) and (1.8) for $\vartheta\in[0,3/2-\alpha]$. ☐

For the remaining case, the endpoint decay established above provides the time integrability required to propagate the negative Sobolev norm.

*Proof of* (1.6)–(1.8) *for* $\vartheta\in(3/2-\alpha,3/2)$. The preceding argument for propagating the negative Sobolev norms does not directly apply when $\vartheta\in(3/2-\alpha,3/2)$. However, since $\dot{H}^{-\vartheta}\cap L^2\subset\dot{H}^{-\vartheta'}$ for any $\vartheta'\in[0,\vartheta]$, the initial data satisfy $\varrho_0,u_0,\Lambda^\beta\phi_0\in\dot{H}^{-(3/2-\alpha)}$. Thus, the decay result (1.7) established for $\vartheta=3/2-\alpha$ implies for $l=0,\dots,N-1$ that

$$\sum_{k=l}^{N}\left(\|\Lambda^{k\alpha}\varrho(t)\|_{L^2}+\|\Lambda^{k\alpha}u(t)\|_{L^2}+\|\Lambda^{k\alpha+\beta}\phi(t)\|_{L^2}\right) \leq C\,(M_0+\delta_0)\,(1+t)^{-\left(\frac{l-1}{2}+\frac{3}{4\alpha}\right)}. \tag{4.9}$$

Consequently, utilizing (3.3) and (4.9), we obtain

$$\|(\varrho,u,\Lambda^\beta\phi)\|_{L^2}^{1-\frac{3-2\vartheta}{2\alpha}}\left(\|\Lambda^\alpha(\varrho,u,\Lambda^\beta\phi)\|_{L^2}+\|\Lambda^{2\alpha}(\varrho,u)\|_{L^2}\right)^{1+\frac{3-2\vartheta}{2\alpha}} \leq C\,(M_0+\delta_0)^2\,(1+t)^{-\frac{9-2\alpha-2\vartheta}{4\alpha}},$$

which is integrable in time since $(9-2\alpha-2\vartheta)/(4\alpha)>1$. Hence, integrating (3.1), we deduce that for $\vartheta\in(3/2-\alpha,3/2)$,

$$\mathcal{E}_{-\vartheta}(t) \leq CM_0^2 + C\,(M_0+\delta_0)^2\sup_{0\leq\tau\leq t}\sqrt{\mathcal{E}_{-\vartheta}(\tau)}.$$

This implies (1.6) for $\vartheta\in(3/2-\alpha,3/2)$

$$\|\Lambda^{-\vartheta}\varrho(t)\|_{L^2}+\|\Lambda^{-\vartheta}u(t)\|_{L^2}+\|\Lambda^{-\vartheta+\beta}\phi(t)\|_{L^2} \leq C\,(M_0+\delta_0)\,(1+M_0+\delta_0)\,. \tag{4.10}$$

Combining (4.10) with the interpolation argument used in the preceding case, we obtain (1.7) and (1.8) for $\vartheta\in(3/2-\alpha,3/2)$. ☐

## 5. Optimality of the decay rates

In this section, we first prove Theorem 1.6 by constructing special initial data satisfying (1.3) and (1.5) and establishing two-sided decay estimates for the corresponding solutions to (1.1). We then use these solutions to prove the optimality assertion in Theorem 1.1(iii). The construction is based on the auxiliary variables

$$w := \Lambda^{-\beta}\varrho = \Lambda^{\beta}\phi, \qquad d := \Lambda^{-1}\mathrm{div}u, \qquad v := \mathbb{P}_{\mathrm{T}}u := (\mathbb{I} - \nabla\Delta^{-1}\mathrm{div})u. \tag{5.1}$$

Here $\mathbb{P}_{\mathrm{T}}$ is the transverse Helmholtz (Leray) projector and $\mathbb{P}_{\mathrm{L}} = \mathbb{I} - \mathbb{P}_{\mathrm{T}}$ is the longitudinal projector. Conversely,

$$\varrho = \Lambda^{\beta}w, \qquad u = -\nabla\Lambda^{-1}d + v, \qquad \Lambda^{\beta}\phi = w. \tag{5.2}$$

Using the variables defined in (5.1), we rewrite the system (2.1) as

$$\begin{cases} \partial_t w + \Lambda^{1-\beta}d = -\Lambda^{-\beta}\mathrm{div}\mathcal{M}, \\ \partial_t d + \nu_{\mathrm{L}}\Lambda^{2\alpha}d - \Lambda^{1-\beta}(\mathbb{I} + \Lambda^{2\beta})w = \Lambda^{-1}\mathrm{div}\mathcal{G}, \\ \partial_t v + \nu_{\mathrm{T}}\Lambda^{2\alpha}v = \mathbb{P}_{\mathrm{T}}\mathcal{G}, \end{cases} \tag{5.3}$$

where the nonlinear source terms are

$$\mathcal{M} := \varrho u, \qquad \mathcal{G} := -u\cdot\nabla u + H(\varrho)\mathfrak{L}^{\alpha}u - F(\varrho)\nabla\varrho.$$

For $r > 0$, define the matrix

$$A(r) := \begin{pmatrix} 0 & r^{1-\beta} \\ -r^{1-\beta}(1 + r^{2\beta}) & \nu_{\mathrm{L}}r^{2\alpha} \end{pmatrix}.$$

Let $\mathcal{S}(t)$ denote the homogeneous linear solution family defined by the following Fourier multipliers. For $U_0 = (w_0, d_0, v_0)$, its component solution operators $\mathcal{S}_w(t)$, $\mathcal{S}_d(t)$, and $\mathcal{S}_v(t)$ are given, for $\xi \neq 0$ with $r = |\xi|$, by

$$\begin{pmatrix} \widehat{\mathcal{S}_w(t)U_0}(\xi) \\ \widehat{\mathcal{S}_d(t)U_0}(\xi) \end{pmatrix} := e^{-tA(r)} \begin{pmatrix} \widehat{w}_0(\xi) \\ \widehat{d}_0(\xi) \end{pmatrix}, \qquad \widehat{\mathcal{S}_v(t)U_0}(\xi) := e^{-\nu_{\mathrm{T}}tr^{2\alpha}}\widehat{v}_0(\xi).$$

At $\xi = 0$, we take the multiplier defining $\mathcal{S}(t)$ to be the identity. Thus

$$\mathcal{S}(t)U_0 = (\mathcal{S}_w(t)U_0,\, \mathcal{S}_d(t)U_0,\, \mathcal{S}_v(t)U_0).$$

For a solution of (5.3) with initial data $U_0$, write $U = (w, d, v)$ and define the remainder

$$R(t) := U(t) - \mathcal{S}(t)U_0 = (R_w(t), R_d(t), R_v(t)).$$

### 5.1. Spectral properties and construction of special initial data.

We first collect the low-frequency spectral and Green matrix results used in the construction of special initial data and the Duhamel analysis.

**Lemma 5.1.** *Let $\alpha \in [1/2, 1]$ and $\beta \in [0, 1]$. There exist constants $0 < c \le C$ such that, for any $0 < r \le 1$, the following assertions hold.*

(i) *The two eigenvalues of $-A(r)$, counted with multiplicity, are*

$$\mu_{\pm}(r) = -\frac{\nu_{\mathrm{L}}}{2}r^{2\alpha} \pm \frac{1}{2}\sqrt{\nu_{\mathrm{L}}^2 r^{4\alpha} - 4r^{2-2\beta}(1 + r^{2\beta})}, \tag{5.4}$$

*where the square root is taken in the principal complex branch. They satisfy the uniform estimates*

$$\begin{aligned} cr^{2\alpha} &\le -\mathrm{Re}\mu_{\pm}(r) \le Cr^{2\alpha}, \\ cr^{1-\beta} &\le |\mu_{\pm}(r)| \le Cr^{1-\beta}. \end{aligned} \tag{5.5}$$

(ii) *Define*

$$E_\pm(r) = \begin{pmatrix} r^{-\beta} \\ -\mu_\pm(r)/r \end{pmatrix}.$$

*These are eigenvectors of $-A(r)$ associated with $\mu_\pm(r)$, and*

$$e^{-tA(r)}E_\pm(r) = e^{\mu_\pm(r)t}E_\pm(r), \qquad t \geq 0. \tag{5.6}$$

*This identity also holds at a repeated root.*

(iii) *For $t \geq 0$, let $g_{ij}(r,t)$ denote the entries of $e^{-tA(r)}$. If $\mu_+(r) \neq \mu_-(r)$, then*

$$\begin{aligned} g_{11}(r,t) &= \frac{\mu_+(r)e^{\mu_-(r)t} - \mu_-(r)e^{\mu_+(r)t}}{\mu_+(r) - \mu_-(r)}, \\ g_{12}(r,t) &= -r^{1-\beta}\frac{e^{\mu_+(r)t} - e^{\mu_-(r)t}}{\mu_+(r) - \mu_-(r)}, \\ g_{21}(r,t) &= r^{1-\beta}(1+r^{2\beta})\frac{e^{\mu_+(r)t} - e^{\mu_-(r)t}}{\mu_+(r) - \mu_-(r)}, \\ g_{22}(r,t) &= \frac{\mu_+(r)e^{\mu_+(r)t} - \mu_-(r)e^{\mu_-(r)t}}{\mu_+(r) - \mu_-(r)}. \end{aligned} \tag{5.7}$$

*At a repeated root, $\mu_+(r) = \mu_-(r) = -\nu_{\mathrm{L}} r^{2\alpha}/2$, and*

$$(g_{ij}(r,t))_{i,j=1}^2 = e^{-\nu_{\mathrm{L}} r^{2\alpha} t/2}\begin{pmatrix} 1 + \nu_{\mathrm{L}} r^{2\alpha}t/2 & -r^{1-\beta}t \\ r^{1-\beta}(1+r^{2\beta})t & 1 - \nu_{\mathrm{L}} r^{2\alpha}t/2 \end{pmatrix}. \tag{5.8}$$

*In all cases,*

$$|g_{ij}(r,t)| \leq Ce^{-cr^{2\alpha}t}, \qquad i,j = 1,2. \tag{5.9}$$

*Proof.* The explicit formulas and the eigenvector identity follow by direct calculation. We first prove (5.5). If $\nu_{\mathrm{L}}^2 r^{4\alpha} \leq 4r^{2-2\beta}(1+r^{2\beta})$, then

$$-\mathrm{Re}\mu_\pm(r) = \frac{\nu_{\mathrm{L}}}{2}r^{2\alpha}, \qquad |\mu_\pm(r)| = r^{1-\beta}\sqrt{1+r^{2\beta}}.$$

Otherwise, both roots are real and $\mu_-(r) < \mu_+(r) < 0$. Vieta's formulas yield

$$\frac{r^{2-2\beta}(1+r^{2\beta})}{\nu_{\mathrm{L}} r^{2\alpha}} \leq -\mu_+(r) \leq -\mu_-(r) \leq \nu_{\mathrm{L}} r^{2\alpha}.$$

Since $2\alpha + \beta - 1 \geq 0$ and $0 < r \leq 1$, these relations prove both estimates in (5.5).

It remains to prove (5.9). We consider the following four cases.

**Case 1:** $\nu_{\mathrm{L}}^2 r^{4\alpha} \leq 2r^{2-2\beta}(1+r^{2\beta})$. The roots are complex conjugates, with $-\mathrm{Re}\mu_\pm(r) = \nu_{\mathrm{L}} r^{2\alpha}/2$. The formulas (5.7) and the triangle inequality give

$$\begin{aligned} |g_{jj}(r,t)| &\leq \frac{|\mu_+(r)| + |\mu_-(r)|}{|\mu_+(r) - \mu_-(r)|}e^{-\nu_{\mathrm{L}} r^{2\alpha}t/2}, \qquad j = 1,2, \\ |g_{12}(r,t)| &\leq \frac{2r^{1-\beta}}{|\mu_+(r) - \mu_-(r)|}e^{-\nu_{\mathrm{L}} r^{2\alpha}t/2}, \\ |g_{21}(r,t)| &\leq \frac{2r^{1-\beta}(1+r^{2\beta})}{|\mu_+(r) - \mu_-(r)|}e^{-\nu_{\mathrm{L}} r^{2\alpha}t/2}. \end{aligned}$$

Moreover,

$$\begin{aligned} |\mu_\pm(r)| &= r^{1-\beta}\sqrt{1+r^{2\beta}}, \\ |\mu_+(r) - \mu_-(r)| &\geq \sqrt{2}\, r^{1-\beta}\sqrt{1+r^{2\beta}}. \end{aligned}$$

Together with $1 \leq 1 + r^{2\beta} \leq 2$, these relations imply (5.9) in this case.

**Case 2:** $2r^{2-2\beta}(1+r^{2\beta}) < \nu_{\mathrm{L}}^2 r^{4\alpha} < 4r^{2-2\beta}(1+r^{2\beta})$. Rewriting the diagonal numerators and using $|\mu_\pm(r)| = r^{1-\beta}\sqrt{1+r^{2\beta}}$ together with $1 \le 1+r^{2\beta} \le 2$, we obtain from (5.7) and the triangle inequality

$$
\begin{aligned}
|g_{11}(r,t)| &\le |e^{\mu_+(r)t}| + |\mu_+(r)| \left| \frac{e^{\mu_+(r)t} - e^{\mu_-(r)t}}{\mu_+(r) - \mu_-(r)} \right|, \\
|g_{22}(r,t)| &\le |e^{\mu_+(r)t}| + |\mu_-(r)| \left| \frac{e^{\mu_+(r)t} - e^{\mu_-(r)t}}{\mu_+(r) - \mu_-(r)} \right|, \\
|g_{12}(r,t)| &\le |\mu_+(r)| \left| \frac{e^{\mu_+(r)t} - e^{\mu_-(r)t}}{\mu_+(r) - \mu_-(r)} \right|, \\
|g_{21}(r,t)| &\le \sqrt{2}\, |\mu_+(r)| \left| \frac{e^{\mu_+(r)t} - e^{\mu_-(r)t}}{\mu_+(r) - \mu_-(r)} \right|.
\end{aligned}
$$

By (5.4) and the left-hand inequality defining this case,

$$
-\mathrm{Re}\mu_\pm(r) = \frac{\nu_{\mathrm{L}}}{2} r^{2\alpha}, \qquad |\mu_\pm(r)| < \frac{\nu_{\mathrm{L}}}{\sqrt{2}} r^{2\alpha}.
$$

Moreover, the sine representation gives

$$
\left| \frac{e^{\mu_+(r)t} - e^{\mu_-(r)t}}{\mu_+(r) - \mu_-(r)} \right| = e^{-\nu_{\mathrm{L}} r^{2\alpha} t/2} \left| \frac{\sin(t\, \mathrm{Im}\, \mu_+(r))}{\mathrm{Im}\, \mu_+(r)} \right| \le t e^{-\nu_{\mathrm{L}} r^{2\alpha} t/2}.
$$

Substituting these bounds into the preceding estimates and absorbing the resulting linear factor into the exponential yields (5.9) in this case.
**Case 3:** $\nu_{\mathrm{L}}^2 r^{4\alpha} = 4r^{2-2\beta}(1+r^{2\beta})$. By (5.8), the entries satisfy

$$
\begin{aligned}
|g_{jj}(r,t)| &\le |e^{\mu_+(r)t}| + |\mu_+(r)|\, |te^{\mu_+(r)t}|, \qquad j = 1,2, \\
|g_{12}(r,t)| &\le |\mu_+(r)|\, |te^{\mu_+(r)t}|, \\
|g_{21}(r,t)| &\le \sqrt{2}\, |\mu_+(r)|\, |te^{\mu_+(r)t}|.
\end{aligned}
$$

Here,

$$
\mu_+(r) = \mu_-(r) = -\frac{\nu_{\mathrm{L}}}{2} r^{2\alpha}, \qquad |te^{\mu_+(r)t}| = te^{-\nu_{\mathrm{L}} r^{2\alpha} t/2}.
$$

Substituting these identities into the preceding estimates and absorbing the resulting linear factor into the exponential yields (5.9) in this case.
**Case 4:** $\nu_{\mathrm{L}}^2 r^{4\alpha} > 4r^{2-2\beta}(1+r^{2\beta})$. Both roots are real and $\mu_-(r) < \mu_+(r) < 0$. Rewriting the diagonal numerators and applying the triangle inequality to (5.7), we obtain

$$
\begin{aligned}
|g_{jj}(r,t)| &\le |e^{\mu_+(r)t}| + |\mu_-(r)| \left| \frac{e^{\mu_+(r)t} - e^{\mu_-(r)t}}{\mu_+(r) - \mu_-(r)} \right|, \qquad j = 1,2, \\
|g_{12}(r,t)| &\le |\mu_-(r)| \left| \frac{e^{\mu_+(r)t} - e^{\mu_-(r)t}}{\mu_+(r) - \mu_-(r)} \right|, \\
|g_{21}(r,t)| &\le \sqrt{2}\, |\mu_-(r)| \left| \frac{e^{\mu_+(r)t} - e^{\mu_-(r)t}}{\mu_+(r) - \mu_-(r)} \right|.
\end{aligned}
$$

By the real-root estimates established above, we have

$$
-\mathrm{Re}\mu_-(r) \ge -\mathrm{Re}\mu_+(r) \ge \frac{r^{2\alpha}}{\nu_{\mathrm{L}}}, \qquad |\mu_\pm(r)| \le \nu_{\mathrm{L}} r^{2\alpha}.
$$

Moreover, the mean value theorem gives

$$\left|\frac{e^{\mu_+(r)t}-e^{\mu_-(r)t}}{\mu_+(r)-\mu_-(r)}\right| \le te^{\mu_+(r)t} \le te^{-r^{2\alpha}t/\nu_{\mathrm{L}}}.$$

Substituting these bounds into the preceding estimates and absorbing the resulting linear factor into the exponential completes the proof of (5.9). □

We next construct real-valued initial data. Fix $r_* \in (0,1]$, and choose $\chi_* \in C_c^\infty([0,\infty))$ with $0 \le \chi_* \le 1$, support in $[0, r_*/2)$, and $\chi_* = 1$ on $[0, r_*/4]$. Write $\xi = r\omega$, where $r = |\xi| > 0$ and $\omega \in \mathbb{S}^2$. Choose a nonzero real angular cutoff function $\eta \in C^\infty(\mathbb{S}^2)$ supported in $\{\omega_1 > 1/2\}$, so that $\eta(\omega)$ and $\eta(-\omega)$ have disjoint supports. Fix $0 \le \vartheta < 3/2$ and $0 < \varepsilon \le 3/2 - \vartheta$. We pair $E_+(r)$ with its complex conjugate to preserve the Fourier symmetry of real-valued data. For $0 < r < r_*/2$, define

$$\begin{pmatrix}\widehat{w}_0\\ \widehat{d}_0\end{pmatrix} := \delta\chi_*(r)r^{\beta+\vartheta+\varepsilon-3/2}\left[\eta(\omega)E_+(r) + \eta(-\omega)\overline{E_+(r)}\right], \qquad \widehat{v}_0 = 0, \tag{5.10}$$

and set the transforms to zero elsewhere and at $\xi = 0$. The corresponding initial data in the original variables are then given by

$$\begin{aligned}\widehat{\varrho}_0 &= \delta\chi_*(r)r^{\beta+\vartheta+\varepsilon-3/2}[\eta(\omega)+\eta(-\omega)],\\ \widehat{u}_0 &= \mathrm{i}\delta\chi_*(r)r^{\beta+\vartheta+\varepsilon-3/2}\frac{\xi}{r^2}[\mu_+(r)\eta(\omega) + \overline{\mu_+(r)}\eta(-\omega)],\\ \widehat{\Lambda^\beta\phi_0} &= r^{-\beta}\widehat{\varrho}_0.\end{aligned} \tag{5.11}$$

**Remark 5.2.** *The construction in* (5.10) *also covers repeated roots, since it uses* $E_+(r)$ *and its complex conjugate and only requires the semigroup identity* (5.6).

The following lemma verifies that these initial data are admissible. More precisely, it establishes the required positive-order Sobolev smallness and negative-order Sobolev regularity. Fix $\sigma \ge 0$ and an integer $N_*$ satisfying (1.9).

**Lemma 5.3.** *Let* $0 \le \vartheta < 3/2$ *and* $0 < \varepsilon \le 3/2 - \vartheta$. *Then there exists* $0 < \delta_2 \le 1$ *such that, for all* $0 < \delta \le \delta_2$, *the initial data constructed by* (5.11) *satisfy*

$$\sum_{k=0}^{N_*}\left(\|\Lambda^{k\alpha}(\rho_0^{\delta,\varepsilon}-1)\|_{L^2}^2 + \|\Lambda^{k\alpha}u_0^{\delta,\varepsilon}\|_{L^2}^2 + \|\Lambda^{k\alpha+\beta}\phi_0^{\delta,\varepsilon}\|_{L^2}^2\right) \le C\delta^2 \le \delta_0^2. \tag{5.12}$$

*Moreover, for any* $0 \le \vartheta' < \vartheta + \varepsilon$,

$$\|\Lambda^{-\vartheta'}(\rho_0^{\delta,\varepsilon}-1)\|_{L^2}^2 + \|\Lambda^{-\vartheta'}u_0^{\delta,\varepsilon}\|_{L^2}^2 + \|\Lambda^{-\vartheta'+\beta}\phi_0^{\delta,\varepsilon}\|_{L^2}^2 \le C\delta^2. \tag{5.13}$$

*Proof.* Since the two angular supports are disjoint, (5.11) and (5.5) give

$$\begin{aligned}|\widehat{\varrho}_0|^2 &= \delta^2\chi_*(r)^2r^{2\beta+2\vartheta+2\varepsilon-3}\left[\eta(\omega)^2+\eta(-\omega)^2\right],\\ |\widehat{\Lambda^\beta\phi_0}|^2 &= \delta^2\chi_*(r)^2r^{2\vartheta+2\varepsilon-3}\left[\eta(\omega)^2+\eta(-\omega)^2\right],\\ c^2\delta^2\chi_*(r)^2r^{2\vartheta+2\varepsilon-3}\left[\eta(\omega)^2+\eta(-\omega)^2\right] &\le |\widehat{u}_0|^2 \le C^2\delta^2\chi_*(r)^2r^{2\vartheta+2\varepsilon-3}\left[\eta(\omega)^2+\eta(-\omega)^2\right].\end{aligned} \tag{5.14}$$

Consequently, for every finite $m \ge 0$,

$$\|\varrho_0^{\delta,\varepsilon}\|_{H^m}^2 + \|u_0^{\delta,\varepsilon}\|_{H^m}^2 + \|\Lambda^\beta\phi_0^{\delta,\varepsilon}\|_{H^m}^2 \le C\delta^2\int_0^{r_*/2}\left(r^{2\beta+2\vartheta+2\varepsilon-1} + r^{2\vartheta+2\varepsilon-1}\right)\mathrm{d}r \le C\delta^2.$$

Taking $m = N_*\alpha$ in the preceding estimate and reducing $\delta_2$ if necessary proves the first estimate. Since $N_*\alpha > 5/2 > 3/2$, Sobolev embedding gives $\|\varrho_0\|_{L^\infty} \le C\delta$. Choosing $\delta_2$ sufficiently small also ensures $\rho_0 \ge 1/2$. The compatibility relation $\Lambda^{2\beta}\phi_0 = \varrho_0$ holds by construction.

Moreover, for any $0 \le \vartheta' < \vartheta + \varepsilon$,

$$\|\varrho_0^{\delta,\varepsilon}\|^2_{\dot H^{-\vartheta'}} + \|u_0^{\delta,\varepsilon}\|^2_{\dot H^{-\vartheta'}} + \|\Lambda^\beta \phi_0^{\delta,\varepsilon}\|^2_{\dot H^{-\vartheta'}} \le C\delta^2 \int_0^{r_*/2} \left(r^{2\beta+2(\vartheta+\varepsilon-\vartheta')-1} + r^{2(\vartheta+\varepsilon-\vartheta')-1}\right) \mathrm{d}r \le C\delta^2.$$

This proves (5.13). ☐

**Remark 5.4.** *The estimate* (5.12) *gives* (1.3) *with $N = N_*$, while* (5.13) *with $\vartheta' = \vartheta$ shows that the left-hand side of* (1.5) *is bounded by $C\delta$, because $\varepsilon > 0$. Thus the data satisfy the hypotheses of Theorem 1.1*(i)–(ii)*.*

**Remark 5.5.** *The range $0 \le \vartheta' < \vartheta + \varepsilon$ in Lemma 5.3 is sharp. At the borderline negative index $\vartheta' = \vartheta + \varepsilon$, $u_0^{\delta,\varepsilon} \notin \dot H^{-(\vartheta+\varepsilon)}$ and $\Lambda^\beta \phi_0^{\delta,\varepsilon} \notin \dot H^{-(\vartheta+\varepsilon)}$. At this endpoint, the weighted Fourier integrals defining these norms diverge logarithmically near $\xi = 0$.*

**Remark 5.6.** *The special data in Lemma 5.3 do not satisfy the $L^1$ assumption imposed in Theorem 1.1 of Wang and Zhang [37]. More precisely, the initial momentum $m_0^{\delta,\varepsilon} := \rho_0^{\delta,\varepsilon} u_0^{\delta,\varepsilon}$ does not belong to $L^1(\mathbb{R}^3)$. Separately, the density perturbation $\rho_0^{\delta,\varepsilon} - 1 = \varrho_0^{\delta,\varepsilon}$ does not belong to $L^1(\mathbb{R}^3)$ whenever $\beta + \vartheta + \varepsilon \le 3/2$. These facts follow directly from the low-frequency profiles in* (5.14)*, which are incompatible with the boundedness and continuity at $\xi = 0$ of the Fourier transform of an $L^1$ function.*

5.2. **Nonlinear decay with explicit amplitude dependence.** Let $(\varrho, u, \Lambda^\beta \phi)$ be the nonlinear solution to (2.1) with initial data given by (5.11). For simplicity, throughout the remainder of this section we omit the superscripts $(\delta, \varepsilon)$ from the nonlinear solution and its sources.

**Lemma 5.7.** *Under the assumptions of Lemma 5.3, let $\delta_2$ be the constant specified therein and choose $\vartheta_0$ such that $(\vartheta + \varepsilon)/2 < \vartheta_0 < \vartheta + \varepsilon$. Then, for any $0 < \delta \le \delta_2$, the following estimates hold:*

(i) *For any $0 \le \sigma' \le (N_* - 1)\alpha$, we have*

$$\|\Lambda^{\sigma'} u(t)\|_{L^2} + \|\Lambda^{\sigma'+\beta} \phi(t)\|_{L^2} \le C\delta(1+t)^{-\frac{\vartheta_0+\sigma'}{2\alpha}}. \tag{5.15}$$

(ii) *For any $0 \le \sigma' \le (N_* - 1)\alpha - \beta$, we have*

$$\|\Lambda^{\sigma'} \varrho(t)\|_{L^2} \le C\delta(1+t)^{-\frac{\vartheta_0+\beta+\sigma'}{2\alpha}}. \tag{5.16}$$

(iii) *For any $0 \le \sigma' \le (N_* - 1)\alpha - 2$, we have*

$$\|\Lambda^{\sigma'} u(t)\|_{L^\infty} \le C\delta(1+t)^{-\frac{\vartheta_0+\sigma'+3/2}{2\alpha}}. \tag{5.17}$$

(iv) *For any $0 \le \sigma' \le (N_* - 1)\alpha - 2 - \beta$, we have*

$$\|\Lambda^{\sigma'} \varrho(t)\|_{L^\infty} \le C\delta(1+t)^{-\frac{\vartheta_0+\beta+\sigma'+3/2}{2\alpha}}. \tag{5.18}$$

*Proof.* By Lemma 5.3, the constructed initial data satisfy the assumptions of Theorem 1.1 with $N = N_*$ and $\vartheta = \vartheta_0$. Moreover, it follows from (5.12) and (5.13) that the corresponding positive-order and negative Sobolev norms of the initial data are bounded by $C\delta$. Therefore, applying the same decay argument as in the proof of Theorem 1.1, we obtain

$$\|\Lambda^{l\alpha} u(t)\|_{L^2} + \|\Lambda^{l\alpha+\beta} \phi(t)\|_{L^2} \le C\delta(1+t)^{-\frac{\vartheta_0+l\alpha}{2\alpha}}, \quad \text{for } l = 0, \ldots, N_* - 1. \tag{5.19}$$

Let $0 \le \sigma' \le (N_* - 1)\alpha$ and set

$$\theta = \frac{\sigma'}{(N_* - 1)\alpha} \in [0, 1].$$

By the interpolation inequality and (5.19), we have

$$\begin{aligned}\|\Lambda^{\sigma'}u(t)\|_{L^2} &\le \|u(t)\|_{L^2}^{1-\theta}\|\Lambda^{(N_*-1)\alpha}u(t)\|_{L^2}^{\theta}\\ &\le C\delta(1+t)^{-\frac{(1-\theta)\vartheta_0+\theta\{\vartheta_0+(N_*-1)\alpha\}}{2\alpha}}\\ &= C\delta(1+t)^{-\frac{\vartheta_0+\sigma'}{2\alpha}}.\end{aligned}$$

Applying the same interpolation to $\Lambda^\beta\phi$ gives

$$\|\Lambda^{\sigma'+\beta}\phi(t)\|_{L^2} \le C\delta(1+t)^{-\frac{\vartheta_0+\sigma'}{2\alpha}}.$$

Since $\varrho = \Lambda^{2\beta}\phi$, for $0 \le \sigma' \le (N_*-1)\alpha - \beta$, we apply (5.15) at order $\sigma'+\beta$ to obtain

$$\|\Lambda^{\sigma'}\varrho(t)\|_{L^2} = \|\Lambda^{\sigma'+2\beta}\phi(t)\|_{L^2} \le C\delta(1+t)^{-\frac{\vartheta_0+\beta+\sigma'}{2\alpha}}.$$

Finally, for $0 \le \sigma' \le (N_*-1)\alpha - 2$, we use the Gagliardo–Nirenberg inequality and (5.15) to deduce

$$\|\Lambda^{\sigma'}u(t)\|_{L^\infty} \le C\|\Lambda^{\sigma'+1}u(t)\|_{L^2}^{1/2}\|\Lambda^{\sigma'+2}u(t)\|_{L^2}^{1/2} \le C\delta(1+t)^{-\frac{\vartheta_0+\sigma'+3/2}{2\alpha}},$$

which proves (5.17). Similarly, (5.18) follows from the same Gagliardo–Nirenberg inequality and (5.16) for $0 \le \sigma' \le (N_*-1)\alpha - 2 - \beta$. This completes the proof. ☐

5.3. **Time-decay estimates for the auxiliary variables.** Write the source vector in (5.3) as

$$\mathcal{B}(\mathcal{M},\mathcal{G}) = (-\Lambda^{-\beta}\mathrm{div}\mathcal{M}, \Lambda^{-1}\mathrm{div}\mathcal{G}, \mathbb{P}_{\mathrm{T}}\mathcal{G}).$$

Applying Duhamel's principle to (5.3), we obtain

$$U(t) = \mathcal{S}(t)U_0 + \int_0^t \mathcal{S}(t-\tau)\mathcal{B}(\mathcal{M},\mathcal{G})(\tau)\,\mathrm{d}\tau. \tag{5.20}$$

By (5.6) and the reality of $A(r)$, the special linear orbit satisfies, for $0 < r < r_*/2$,

$$\begin{pmatrix}\widehat{\mathcal{S}_w(t)U_0}\\ \widehat{\mathcal{S}_d(t)U_0}\end{pmatrix} = \delta\chi_*(r)r^{\beta+\vartheta+\varepsilon-3/2}\big[e^{\mu_+(r)t}\eta(\omega)E_+(r) + e^{\overline{\mu_+(r)}t}\eta(-\omega)\overline{E_+(r)}\big], \qquad \widehat{\mathcal{S}_v(t)U_0} = 0. \tag{5.21}$$

**Lemma 5.8.** *For any $\sigma' \ge 0$, there are constants $0 < c \le C$ such that $\mathcal{S}_v(t)U_0 = 0$ and*

$$c\delta(1+t)^{-\frac{\sigma'+\vartheta+\varepsilon}{2\alpha}} \le \|\Lambda^{\sigma'}\mathcal{S}_w(t)U_0\|_{L^2} \le C\delta(1+t)^{-\frac{\sigma'+\vartheta+\varepsilon}{2\alpha}}, \tag{5.22}$$

$$c\delta(1+t)^{-\frac{\sigma'+\vartheta+\varepsilon}{2\alpha}} \le \|\Lambda^{\sigma'}\mathcal{S}_d(t)U_0\|_{L^2} \le C\delta(1+t)^{-\frac{\sigma'+\vartheta+\varepsilon}{2\alpha}}. \tag{5.23}$$

*Proof.* Since the two angular caps are disjoint and complex conjugation preserves absolute values, (5.21) gives

$$\begin{aligned}|\widehat{\mathcal{S}_w(t)U_0}|^2 &= \delta^2\chi_*^2 r^{2(\vartheta+\varepsilon)-3}e^{2\mathrm{Re}\mu_+(r)t}\big[\eta(\omega)^2+\eta(-\omega)^2\big],\\ |\widehat{\mathcal{S}_d(t)U_0}|^2 &= \delta^2\chi_*^2 r^{2\beta+2(\vartheta+\varepsilon)-5}|\mu_+(r)|^2e^{2\mathrm{Re}\mu_+(r)t}\big[\eta(\omega)^2+\eta(-\omega)^2\big].\end{aligned} \tag{5.24}$$

Set $\zeta = \sigma'+\vartheta+\varepsilon > 0$. For $a > 0$, define

$$I_a(t) := \int_0^{r_*/2}\chi_*(r)^2 r^{2\zeta-1}e^{-atr^{2\alpha}}\,\mathrm{d}r.$$

Applying Plancherel's theorem to (5.24) and using (5.5), we obtain, for suitable fixed constants $0 < c \le C$,

$$c\delta^2 I_{2C}(t) \le \|\Lambda^{\sigma'}\mathcal{S}_w(t)U_0\|_{L^2}^2 \le C\delta^2 I_{2c}(t), \tag{5.25}$$

$$c\delta^2 I_{2C}(t) \le \|\Lambda^{\sigma'}\mathcal{S}_d(t)U_0\|_{L^2}^2 \le C\delta^2 I_{2c}(t). \tag{5.26}$$

It therefore suffices to estimate $I_a(t)$ for fixed $a > 0$. We first establish the upper bound. Since $0 \le \chi_* \le 1$ and $\operatorname{supp}\chi_* \subset [0, r_*/2)$, for any $t \ge 0$, we have

$$I_a(t) \le e^{a(r_*/2)^{2\alpha}} \int_0^\infty r^{2\zeta-1} e^{-a(1+t)r^{2\alpha}}\,\mathrm{d}r = \frac{\Gamma(\zeta/\alpha)}{2\alpha} a^{-\zeta/\alpha} e^{a(r_*/2)^{2\alpha}} (1+t)^{-\zeta/\alpha}.$$

For the lower bound, since $\chi_*(r) = 1$ on $[0, r_*/4]$, for any $t \ge 0$, we have

$$I_a(t) \ge e^{-a(r_*/4)^{2\alpha}} \int_0^{(r_*/4)(1+t)^{-1/(2\alpha)}} r^{2\zeta-1}\,\mathrm{d}r = \frac{e^{-a(r_*/4)^{2\alpha}}}{2\zeta} \left(\frac{r_*}{4}\right)^{2\zeta} (1+t)^{-\zeta/\alpha}.$$

Consequently,

$$c(1+t)^{-\zeta/\alpha} \le I_a(t) \le C(1+t)^{-\zeta/\alpha}. \tag{5.27}$$

The bounds (5.25) and (5.26), together with (5.27) for $a = 2c$ and $a = 2C$, yield (5.22) and (5.23) after taking square roots. ☐

By (5.20), the remainder satisfies

$$R(t) = \int_0^t \mathcal{S}(t-\tau)\mathcal{B}(\mathcal{M},\mathcal{G})(\tau)\,\mathrm{d}\tau. \tag{5.28}$$

Choose a radial cutoff $\chi_{\mathrm{L}} \in C_c^\infty(\mathbb{R}^3)$ such that $0 \le \chi_{\mathrm{L}} \le 1$, $\chi_{\mathrm{L}}(\xi) = 1$ for $|\xi| \le 1/2$, and $\chi_{\mathrm{L}}(\xi) = 0$ for $|\xi| \ge 1$. Since $r_* \le 1$, the special orbit is supported in the region where $\chi_{\mathrm{L}} = 1$, and

$$\chi_{\mathrm{L}}(D)U(t) = \mathcal{S}(t)U_0 + \chi_{\mathrm{L}}(D)R(t).$$

**Lemma 5.9.** *Suppose that $0 \le \vartheta < 3/2$, $0 < \varepsilon \le 3/2 - \vartheta$, and $(\vartheta+\varepsilon)/2 < \vartheta_0 < \vartheta + \varepsilon$. Define*

$$\eta_* := \frac{1}{4\alpha} \min\left\{\frac{3}{2} - \vartheta - \varepsilon,\, 2\vartheta_0 - \vartheta - \varepsilon + \frac{5}{2} - 2\alpha\right\}. \tag{5.29}$$

*Then $\eta_* > 0$ for $\vartheta + \varepsilon < 3/2$, and, for any $0 < \delta \le \delta_2$ and $t \ge 0$,*

$$\sum_{\sharp \in \{w,d,v\}} \|\Lambda^{\sigma'} \chi_{\mathrm{L}}(D) R_\sharp(t)\|_{L^2} \le C\delta^2 (1+t)^{-\frac{\sigma'+\vartheta+\varepsilon}{2\alpha} - \eta_*}, \qquad 0 \le \sigma' \le \sigma + \beta. \tag{5.30}$$

*Proof.* We use (5.9) from Lemma 5.1(iii) to derive a low-frequency Fourier estimate for $\mathcal{S}(\tau)\mathcal{B}(\mathcal{M},\mathcal{G})$. For $\xi \ne 0$, write $r = |\xi|$. By the definitions of $\mathcal{B}$ and $\mathcal{S}(\tau)$,

$$\begin{aligned}
\widehat{\mathcal{S}_w(\tau)\mathcal{B}(\mathcal{M},\mathcal{G})}(\xi) &= -\mathrm{i}r^{-\beta} g_{11}(r,\tau)\,\xi \cdot \widehat{\mathcal{M}}(\xi) + \mathrm{i}r^{-1} g_{12}(r,\tau)\,\xi \cdot \widehat{\mathcal{G}}(\xi),\\
\widehat{\mathcal{S}_d(\tau)\mathcal{B}(\mathcal{M},\mathcal{G})}(\xi) &= -\mathrm{i}r^{-\beta} g_{21}(r,\tau)\,\xi \cdot \widehat{\mathcal{M}}(\xi) + \mathrm{i}r^{-1} g_{22}(r,\tau)\,\xi \cdot \widehat{\mathcal{G}}(\xi),\\
\widehat{\mathcal{S}_v(\tau)\mathcal{B}(\mathcal{M},\mathcal{G})}(\xi) &= e^{-\nu_{\mathrm{T}}\tau r^{2\alpha}} \mathbb{P}_{\mathrm{T}}(\xi)\widehat{\mathcal{G}}(\xi).
\end{aligned}$$

Moreover, one has

$$r^{-\beta}|\xi \cdot \widehat{\mathcal{M}}(\xi)| \le r^{1-\beta}|\widehat{\mathcal{M}}(\xi)|, \qquad r^{-1}|\xi \cdot \widehat{\mathcal{G}}(\xi)| \le |\widehat{\mathcal{G}}(\xi)|, \qquad |\mathbb{P}_{\mathrm{T}}(\xi)\widehat{\mathcal{G}}(\xi)| \le |\widehat{\mathcal{G}}(\xi)|.$$

Combining (5.9) with the preceding identities and bounds yields

$$\sum_{\sharp \in \{w,d,v\}} \left|\widehat{\mathcal{S}_\sharp(\tau)\mathcal{B}(\mathcal{M},\mathcal{G})}(\xi)\right| \le Ce^{-cr^{2\alpha}\tau}\left(r^{1-\beta}|\widehat{\mathcal{M}}(\xi)| + |\widehat{\mathcal{G}}(\xi)|\right), \qquad 0 < r \le 1.$$

Fix $0 \le \sigma' \le \sigma + \beta$. Condition (1.9) ensures that Lemma 5.7 applies at all derivative orders used in the following estimates. For $\mathcal{M} = \varrho u$, Hölder's inequality, Lemma A.3 and Lemma 5.7 give, for $s \ge 0$,

$$\begin{aligned}
\|\Lambda^{\sigma'+1-\beta}(\varrho u)(s)\|_{L^2} &\le C(\|\varrho\|_{L^\infty}\|\Lambda^{\sigma'+1-\beta}u\|_{L^2} + \|\Lambda^{\sigma'+1-\beta}\varrho\|_{L^2}\|u\|_{L^\infty})\\
&\le C\delta^2(1+s)^{-\frac{2\vartheta_0+\sigma'+5/2}{2\alpha}}.
\end{aligned}$$

For the fractional Lamé term in $\mathcal{G}$, by the Gagliardo–Nirenberg inequality, Lemma A.1(iii) and Lemma 5.7, we obtain

$$\|\mathbb{B}_\alpha \Lambda^{2\alpha} u(s)\|_{L^\infty} \le C\|\Lambda^{2\alpha} u(s)\|_{L^2}^{1/4} \|\Lambda^{2\alpha+2} u(s)\|_{L^2}^{3/4} \le C\delta(1+s)^{-\frac{\vartheta_0+2\alpha+3/2}{2\alpha}}.$$

Hence, Hölder's inequality, Lemma A.3, (2.6) and (2.7), together with this bound, imply

$$\begin{aligned}\|\Lambda^{\sigma'}[H(\varrho)\mathbb{B}_\alpha \Lambda^{2\alpha} u](s)\|_{L^2} &\le C\big(\|\varrho\|_{L^\infty}\|\Lambda^{\sigma'+2\alpha} u\|_{L^2} + \|\Lambda^{\sigma'}\varrho\|_{L^2}\|\mathbb{B}_\alpha\Lambda^{2\alpha} u\|_{L^\infty}\big)\\ &\le C\delta^2(1+s)^{-\frac{2\vartheta_0+\sigma'+\beta+2\alpha+3/2}{2\alpha}},\end{aligned}$$

$$\begin{aligned}\|\Lambda^{\sigma'}(u\cdot\nabla u)(s)\|_{L^2} &\le C\big(\|u\|_{L^\infty}\|\Lambda^{\sigma'+1}u\|_{L^2} + \|\Lambda^{\sigma'}u\|_{L^2}\|\nabla u\|_{L^\infty}\big)\\ &\le C\delta^2(1+s)^{-\frac{2\vartheta_0+\sigma'+5/2}{2\alpha}},\end{aligned}$$

and

$$\begin{aligned}\|\Lambda^{\sigma'}[F(\varrho)\nabla\varrho](s)\|_{L^2} &\le C\big(\|\varrho\|_{L^\infty}\|\Lambda^{\sigma'+1}\varrho\|_{L^2} + \|\Lambda^{\sigma'}\varrho\|_{L^2}\|\nabla\varrho\|_{L^\infty}\big)\\ &\le C\delta^2(1+s)^{-\frac{2\vartheta_0+\sigma'+2\beta+5/2}{2\alpha}}.\end{aligned}$$

In view of $\beta \ge 0$ and $\beta + 2\alpha \ge 1$, Plancherel's theorem and the preceding Fourier estimate therefore yield, for all $s, \tau \ge 0$,

$$\sum_{\sharp\in\{w,d,v\}} \|\Lambda^{\sigma'}\chi_{\mathrm{L}}(D)\mathcal{S}_\sharp(\tau)\mathcal{B}(\mathcal{M}(s),\mathcal{G}(s))\|_{L^2} \le C\delta^2(1+s)^{-\frac{2\vartheta_0+\sigma'+5/2}{2\alpha}}. \tag{5.31}$$

We next derive the corresponding $L^1$–$L^2$ estimate. Since $\chi_{\mathrm{L}}$ is supported in $\{r \le 1\}$, the preceding pointwise Fourier estimate, Plancherel's theorem, and the Hausdorff–Young inequality $\|\widehat{f}\|_{L^\infty} \le C\|f\|_{L^1}$ give, for any $\sigma' \ge 0$ and $\tau \ge 0$,

$$\begin{aligned}&\sum_{\sharp\in\{w,d,v\}} \|\Lambda^{\sigma'}\chi_{\mathrm{L}}(D)\mathcal{S}_\sharp(\tau)\mathcal{B}(\mathcal{M},\mathcal{G})\|_{L^2}\\ &\le C\left(\int_{|\xi|\le 1} r^{2(\sigma'+1-\beta)}e^{-2cr^{2\alpha}\tau}|\widehat{\mathcal{M}}(\xi)|^2\,\mathrm{d}\xi\right)^{1/2} + C\left(\int_{|\xi|\le 1} r^{2\sigma'}e^{-2cr^{2\alpha}\tau}|\widehat{\mathcal{G}}(\xi)|^2\,\mathrm{d}\xi\right)^{1/2}\\ &\le C\left(\sqrt{J_{\sigma'+1-\beta}(\tau)}\|\mathcal{M}\|_{L^1} + \sqrt{J_{\sigma'}(\tau)}\|\mathcal{G}\|_{L^1}\right),\end{aligned}$$

where

$$J_\gamma(\tau) := \int_{|\xi|\le 1} |\xi|^{2\gamma} e^{-2c|\xi|^{2\alpha}\tau}\,\mathrm{d}\xi, \qquad \gamma \ge 0.$$

Since $e^{-2c\tau r^{2\alpha}} \le e^{2c}e^{-2c(1+\tau)r^{2\alpha}}$ for $0 \le r \le 1$, the change of variables $y = (1+\tau)^{1/(2\alpha)}r$ yields

$$J_\gamma(\tau) \le C(1+\tau)^{-\frac{2\gamma+3}{2\alpha}}\int_0^\infty y^{2\gamma+2}e^{-2cy^{2\alpha}}\,\mathrm{d}y \le C(1+\tau)^{-\frac{2\gamma+3}{2\alpha}}.$$

Taking $\gamma = \sigma' + 1 - \beta$ and $\gamma = \sigma'$ gives

$$\sum_{\sharp\in\{w,d,v\}} \|\Lambda^{\sigma'}\chi_{\mathrm{L}}(D)\mathcal{S}_\sharp(\tau)\mathcal{B}(\mathcal{M},\mathcal{G})\|_{L^2} \le C(1+\tau)^{-\frac{\sigma'+5/2-\beta}{2\alpha}}\|\mathcal{M}\|_{L^1} + C(1+\tau)^{-\frac{\sigma'+3/2}{2\alpha}}\|\mathcal{G}\|_{L^1}. \tag{5.32}$$

To apply (5.32), we estimate the nonlinear source terms in $L^1$. Hölder's inequality, Lemma 5.7, (2.7), and Lemma A.1(iii) yield, for any $s \geq 0$,

$$\begin{aligned}\|\mathcal{M}(s)\|_{L^1} &= \|\varrho u\|_{L^1} \leq \|\varrho\|_{L^2}\|u\|_{L^2} \leq C\delta^2(1+s)^{-\frac{2\vartheta_0+\beta}{2\alpha}},\\ \|\mathcal{G}(s)\|_{L^1} &\leq C\big(\|u\|_{L^2}\|\nabla u\|_{L^2} + \|\varrho\|_{L^2}\|\Lambda^{2\alpha}u\|_{L^2} + \|\varrho\|_{L^2}\|\nabla\varrho\|_{L^2}\big)\\ &\leq C\delta^2\Big[(1+s)^{-\frac{2\vartheta_0+1}{2\alpha}} + (1+s)^{-\frac{2\vartheta_0+\beta+2\alpha}{2\alpha}} + (1+s)^{-\frac{2\vartheta_0+2\beta+1}{2\alpha}}\Big]\\ &\leq C\delta^2(1+s)^{-\frac{2\vartheta_0+1}{2\alpha}},\end{aligned}$$

where in the last inequality we used $\beta + 2\alpha \geq 1$ and $\beta \geq 0$.

We now split the time integral in (5.28) at $t/2$. Applying (5.32) with the above $L^1$ bounds on $[0, t/2]$ and (5.31) on $[t/2, t]$, we obtain

$$\begin{aligned}&\sum_{\sharp\in\{w,d,v\}} \|\Lambda^{\sigma'}\chi_{\mathrm{L}}(D)R_\sharp(t)\|_{L^2}\\ &\leq C\delta^2\Big[\int_0^{t/2}(1+t-\tau)^{-\frac{\sigma'+5/2-\beta}{2\alpha}}(1+\tau)^{-\frac{2\vartheta_0+\beta}{2\alpha}}\,\mathrm{d}\tau + \int_{t/2}^t(1+\tau)^{-\frac{2\vartheta_0+\sigma'+5/2}{2\alpha}}\,\mathrm{d}\tau\Big]\\ &\quad + C\delta^2\Big[\int_0^{t/2}(1+t-\tau)^{-\frac{\sigma'+3/2}{2\alpha}}(1+\tau)^{-\frac{2\vartheta_0+1}{2\alpha}}\,\mathrm{d}\tau + \int_{t/2}^t(1+\tau)^{-\frac{2\vartheta_0+\sigma'+5/2}{2\alpha}}\,\mathrm{d}\tau\Big]. \qquad (5.33)\end{aligned}$$

For the first bracketed term in (5.33), we apply Lemma A.8 with

$$r_1 = \frac{\sigma' + 5/2 - \beta}{2\alpha}, \qquad r_2 = \frac{2\vartheta_0 + \beta}{2\alpha},$$

while for the second bracketed term we take

$$r_1 = \frac{\sigma' + 3/2}{2\alpha}, \qquad r_2 = \frac{2\vartheta_0 + 1}{2\alpha}.$$

For both choices, $r_1, r_2 > 0$ and

$$r_1 + r_2 = \frac{2\vartheta_0 + \sigma' + 5/2}{2\alpha} > 1.$$

Moreover, since $0 \leq \beta \leq 1$, the definition of $\eta_*$ in (5.29) gives

$$\begin{aligned}&\min\Big\{\frac{\sigma' + 5/2 - \beta}{2\alpha}, \frac{\sigma' + 3/2}{2\alpha}, \frac{2\vartheta_0 + \sigma' + 5/2}{2\alpha} - 1\Big\}\\ &= \frac{\sigma' + \vartheta + \varepsilon}{2\alpha} + \frac{1}{2\alpha}\min\Big\{\frac{3}{2} - \vartheta - \varepsilon,\, 2\vartheta_0 - \vartheta - \varepsilon + \frac{5}{2} - 2\alpha\Big\}\\ &= \frac{\sigma' + \vartheta + \varepsilon}{2\alpha} + 2\eta_*.\end{aligned}$$

Set $r = (\sigma' + \vartheta + \varepsilon)/(2\alpha) + \eta_*$. If $\vartheta + \varepsilon < 3/2$, then $\eta_* > 0$, and hence

$$0 < r < \min\{r_1, r_1 + r_2 - 1\}$$

for both choices. Thus, (A.6) applies directly, with the strict gap absorbing any logarithmic factor. If $\vartheta + \varepsilon = 3/2$, then $\eta_* = 0$, and both choices satisfy

$$r \leq r_1, \qquad r < r_1 + r_2 - 1.$$

When $r < r_1$, we again apply (A.6). When $r = r_1$, the second inequality implies $r_2 > 1$, so the endpoint estimate in Lemma A.8 applies.

Consequently, for any $t \geq 0$,

$$\sum_{\sharp\in\{w,d,v\}} \|\Lambda^{\sigma'}\chi_{\mathrm{L}}(D)R_\sharp(t)\|_{L^2} \leq C\delta^2(1+t)^{-\frac{\sigma'+\vartheta+\varepsilon}{2\alpha}-\eta_*}.$$

This proves (5.30). □

### 5.4. **Proof of Theorem 1.6.**

*Proof of Theorem 1.6.* Choose $\vartheta_0$ satisfying $(\vartheta+\varepsilon)/2 < \vartheta_0 < \vartheta+\varepsilon$. Lemma 5.3 verifies the initial-data assumptions of Theorem 1.6 and also gives $H^m$ regularity for every $m \geq 0$. Together with Theorem 1.1(i), it yields a unique global nonlinear solution for $0 < \delta \leq \delta_2$.

Since the linear solution $\mathcal{S}(t)U_0^{\delta,\varepsilon}$ has compact Fourier support in the low-frequency region, the reconstruction formulas (5.2), Plancherel's theorem, and Lemmas 5.8 and 5.9 yield, for sufficiently small $\delta > 0$ and all $t \geq 0$,

$$\begin{aligned}
\|\Lambda^\sigma\varrho(t)\|_{L^2} &\geq \|\Lambda^{\sigma+\beta}\mathcal{S}_w(t)U_0^{\delta,\varepsilon}\|_{L^2} - \|\Lambda^{\sigma+\beta}\chi_{\mathrm{L}}(D)R_w(t)\|_{L^2}\\
&\geq c\delta(1+t)^{-\frac{\sigma+\beta+\vartheta+\varepsilon}{2\alpha}} - C\delta^2(1+t)^{-\frac{\sigma+\beta+\vartheta+\varepsilon}{2\alpha}-\eta_*}\\
&\geq \frac{c}{2}\delta(1+t)^{-\frac{\sigma+\beta+\vartheta+\varepsilon}{2\alpha}},
\end{aligned}$$

and

$$\begin{aligned}
\|\Lambda^\sigma u(t)\|_{L^2} &\geq \|\Lambda^{\sigma}\mathcal{S}_d(t)U_0^{\delta,\varepsilon}\|_{L^2} - \|\Lambda^{\sigma}\chi_{\mathrm{L}}(D)R_d(t)\|_{L^2}\\
&\geq c\delta(1+t)^{-\frac{\sigma+\vartheta+\varepsilon}{2\alpha}} - C\delta^2(1+t)^{-\frac{\sigma+\vartheta+\varepsilon}{2\alpha}-\eta_*}\\
&\geq \frac{c}{2}\delta(1+t)^{-\frac{\sigma+\vartheta+\varepsilon}{2\alpha}}.
\end{aligned}$$

The preceding estimates establish the lower bounds for $\varrho$ and $u$ in Theorem 1.6. The lower bound for $\Lambda^\beta\phi$ follows by applying the same argument to $w$ at order $\sigma$ and using $w = \Lambda^\beta\phi$.

It remains to prove the upper bounds. We first consider the high-frequency part. Choose $\vartheta+\varepsilon-\vartheta_0 < m_{\mathrm{H}} < 1$. By (1.9), $\sigma+m_{\mathrm{H}}+\beta < (N_*-1)\alpha$. Since $1-\chi_{\mathrm{L}}$ is supported in $\{|\xi| \geq 1/2\}$, Plancherel's theorem and (5.15)–(5.16), applied at $\sigma' = \sigma+m_{\mathrm{H}}$, give

$$\|\Lambda^\sigma(1-\chi_{\mathrm{L}})(D)\varrho(t)\|_{L^2} \leq C\|\Lambda^{\sigma+m_{\mathrm{H}}}\varrho(t)\|_{L^2} \leq C\delta(1+t)^{-\frac{\sigma+\beta+\vartheta+\varepsilon}{2\alpha}},$$

and

$$\|\Lambda^\sigma(1-\chi_{\mathrm{L}})(D)u(t)\|_{L^2} \leq C\|\Lambda^{\sigma+m_{\mathrm{H}}}u(t)\|_{L^2} \leq C\delta(1+t)^{-\frac{\sigma+\vartheta+\varepsilon}{2\alpha}}.$$

For the low-frequency part, the reconstruction formulas (5.2) and $\mathcal{S}_v(t)U_0^{\delta,\varepsilon} = 0$, together with Lemma 5.8 and Lemma 5.9, give, for any $t \geq 0$,

$$\begin{aligned}
\|\Lambda^\sigma\chi_{\mathrm{L}}(D)\varrho(t)\|_{L^2} &\leq \|\Lambda^{\sigma+\beta}\mathcal{S}_w(t)U_0^{\delta,\varepsilon}\|_{L^2} + \|\Lambda^{\sigma+\beta}\chi_{\mathrm{L}}(D)R_w(t)\|_{L^2}\\
&\leq C\delta(1+t)^{-\frac{\sigma+\beta+\vartheta+\varepsilon}{2\alpha}} + C\delta^2(1+t)^{-\frac{\sigma+\beta+\vartheta+\varepsilon}{2\alpha}-\eta_*}\\
&\leq C\delta(1+t)^{-\frac{\sigma+\beta+\vartheta+\varepsilon}{2\alpha}},
\end{aligned}$$

and

$$\begin{aligned}
\|\Lambda^\sigma\chi_{\mathrm{L}}(D)u(t)\|_{L^2} &\leq \|\Lambda^{\sigma}\mathcal{S}_d(t)U_0^{\delta,\varepsilon}\|_{L^2} + \sum_{\sharp\in\{d,v\}}\|\Lambda^{\sigma}\chi_{\mathrm{L}}(D)R_\sharp(t)\|_{L^2}\\
&\leq C\delta(1+t)^{-\frac{\sigma+\vartheta+\varepsilon}{2\alpha}} + C\delta^2(1+t)^{-\frac{\sigma+\vartheta+\varepsilon}{2\alpha}-\eta_*}\\
&\leq C\delta(1+t)^{-\frac{\sigma+\vartheta+\varepsilon}{2\alpha}}.
\end{aligned}$$

Combining the high- and low-frequency estimates yields the upper bounds for $\varrho$ and $u$ in Theorem 1.6. The upper bound for $\Lambda^\beta\phi$ follows by applying the same argument to $w$ at order $\sigma$ and using $w = \Lambda^\beta\phi$. All estimates hold for sufficiently small $\delta > 0$ and all $t \geq 0$. □

5.5. **Proof of the optimality assertion in Theorem 1.1.**

*Proof of Theorem 1.1*(iii). Fix $N$ and $\vartheta$ as in Theorem 1.1(iii), $l \in \{0, \dots, N-1\}$, and $\kappa > 0$. Choose

$$0 < \varepsilon < \min\{2\alpha\kappa, 3/2 - \vartheta\}, \qquad \sigma = l\alpha,$$

and take $N_* \geq N$ satisfying (1.9). Apply Theorem 1.6 with these parameters, choosing $\delta > 0$ sufficiently small that the resulting initial data satisfy (1.3) and (1.5) for the fixed $N$. By uniqueness, the resulting solution agrees with that in Theorem 1.1(i). If the velocity decay exponent could be increased by $\kappa$ for all such initial data, then (1.11) with $\sigma = l\alpha$ and the improved upper bound would give

$$c\delta(1+t)^{-\frac{l\alpha+\vartheta+\varepsilon}{2\alpha}} \leq \|\Lambda^{l\alpha}u(t)\|_{L^2} \leq C(1+t)^{-\frac{l\alpha+\vartheta}{2\alpha}-\kappa}, \qquad t \geq 0.$$

Multiplying by $(1+t)^{(l\alpha+\vartheta)/(2\alpha)+\kappa}$, we obtain

$$c\delta(1+t)^{\kappa-\varepsilon/(2\alpha)} \leq C, \qquad t \geq 0.$$

This contradicts $\varepsilon < 2\alpha\kappa$ as $t \to \infty$. The same argument using (1.12) proves the assertion for $\Lambda^\beta\phi$. If $0 \leq l \leq N-2$ and $(N-1-l)\alpha \geq \beta$, then $\min\{\beta, (N-1-l)\alpha\} = \beta$, and (1.10) gives the same contradiction for the density exponent in (1.8). ☐

## Appendix A. Analytic tools

This appendix collects the analytic tools used in the energy and decay estimates. These results provide bounds for the fractional dissipation, nonlinear terms in Sobolev spaces, and time convolutions arising from the Duhamel formula.

We begin with the structural properties of the fractional Lamé–Navier operator.

**Lemma A.1** (Structural properties of the fractional Lamé–Navier operator)**.** *Let $\alpha \in (0,1]$ and set*

$$\nu_- := \min\{\nu_{\mathrm{T}}, \nu_{\mathrm{L}}\}, \qquad \nu_+ := \max\{\nu_{\mathrm{T}}, \nu_{\mathrm{L}}\}.$$

*Then the following assertions hold.*

(i) *The operator $\mathfrak{L}^\alpha$ in* (1.2) *is a nonnegative self-adjoint operator on $L^2(\mathbb{R}^3;\mathbb{R}^3)$ with domain $H^{2\alpha}(\mathbb{R}^3;\mathbb{R}^3)$.*

(ii) *For any $s \in \mathbb{R}$ and $u \in \dot{H}^{s+\alpha}(\mathbb{R}^3;\mathbb{R}^3)$, we have*

$$\nu_-\|\Lambda^{s+\alpha}u\|_{L^2}^2 \leq (\mathbb{B}_\alpha\Lambda^{s+\alpha}u, \Lambda^{s+\alpha}u) \leq \nu_+\|\Lambda^{s+\alpha}u\|_{L^2}^2. \tag{A.1}$$

(iii) *The multiplier $\mathbb{B}_\alpha$ is self-adjoint on $L^2$ and, for any $1 < p < +\infty$, extends to a bounded operator on $L^p(\mathbb{R}^3;\mathbb{R}^3)$. Moreover, for any $\sigma \in \mathbb{R}$,*

$$\mathbb{B}_\alpha\Lambda^\sigma f = \Lambda^\sigma\mathbb{B}_\alpha f$$

*whenever both sides are defined.*

*Proof.* For $\xi \neq 0$, $\mathbb{P}_{\mathrm{L}}(\xi)$ and $\mathbb{P}_{\mathrm{T}}(\xi)$ are mutually orthogonal real symmetric projections, and the Fourier symbol of the classical positive Lamé–Navier operator is

$$|\xi|^2\big(\mu\mathbb{P}_{\mathrm{T}}(\xi) + (2\mu+\lambda)\mathbb{P}_{\mathrm{L}}(\xi)\big).$$

The spectral calculus for these two projections therefore gives

$$\big[|\xi|^2\big(\mu\mathbb{P}_{\mathrm{T}}(\xi) + (2\mu+\lambda)\mathbb{P}_{\mathrm{L}}(\xi)\big)\big]^\alpha = |\xi|^{2\alpha}\big(\mu^\alpha\mathbb{P}_{\mathrm{T}}(\xi) + (2\mu+\lambda)^\alpha\mathbb{P}_{\mathrm{L}}(\xi)\big).$$

This is Scott's spectral representation [28, Theorem 1.1 and Lemma 2.1], written in the Fourier normalization used here. Since its two eigenvalues are $\nu_{\mathrm{T}}|\xi|^{2\alpha}$ and $\nu_{\mathrm{L}}|\xi|^{2\alpha}$, assertion (i) follows from the $L^2$ Fourier functional calculus.

Using Plancherel's theorem and the orthogonality of $\mathbb{P}_{\mathrm{T}}$ and $\mathbb{P}_{\mathrm{L}}$, we obtain

$$(\mathbb{B}_\alpha\Lambda^{s+\alpha}u, \Lambda^{s+\alpha}u) = \nu_{\mathrm{T}}\|\Lambda^{s+\alpha}\mathbb{P}_{\mathrm{T}}u\|_{L^2}^2 + \nu_{\mathrm{L}}\|\Lambda^{s+\alpha}\mathbb{P}_{\mathrm{L}}u\|_{L^2}^2.$$

The bounds in (A.1) follow because $\mathbb{P}_\mathrm{T} + \mathbb{P}_\mathrm{L} = \mathbb{I}$ and $\nu_- \le \nu_\mathrm{T}, \nu_\mathrm{L} \le \nu_+$.

The real symmetry of $\mathbb{P}_\mathrm{T}(\xi)$ and $\mathbb{P}_\mathrm{L}(\xi)$ proves that $\mathbb{B}_\alpha$ is self-adjoint on $L^2$. If $\mathcal{R}_j$ denotes the $j$th Riesz transform, then

$$(\mathbb{P}_\mathrm{L} f)_j = -\sum_{k=1}^{3} \mathcal{R}_j \mathcal{R}_k f_k.$$

The $L^p$-boundedness of the Riesz transforms for $1 < p < +\infty$ implies the corresponding boundedness of $\mathbb{P}_\mathrm{L}$, $\mathbb{P}_\mathrm{T} = \mathbb{I} - \mathbb{P}_\mathrm{L}$, and $\mathbb{B}_\alpha = \nu_\mathrm{T}\mathbb{P}_\mathrm{T} + \nu_\mathrm{L}\mathbb{P}_\mathrm{L}$. Finally, $\mathbb{B}_\alpha$ commutes with $\Lambda^\sigma$ because their Fourier symbols commute pointwise. This proves assertion (iii). □

We next recall the Gagliardo–Nirenberg interpolation inequality used to estimate the nonlinear terms.

**Lemma A.2** ([26])**.** *Let $0 \le s_1, s_2 \le s_3$. For $2 \le p < +\infty$, we have*

$$\|\nabla^{s_1} f\|_{L^p} \le C \|\nabla^{s_2} f\|_{L^2}^{\theta} \|\nabla^{s_3} f\|_{L^2}^{1-\theta},$$

*where $0 \le \theta \le 1$ satisfies*

$$\frac{1}{p} - \frac{s_1}{3} = \left(\frac{1}{2} - \frac{s_2}{3}\right)\theta + \left(\frac{1}{2} - \frac{s_3}{3}\right)(1-\theta).$$

*For $p = +\infty$, if $s_2 \le s_1 + 1$ and $s_3 \ge s_1 + 2$, then*

$$\|\nabla^{s_1} f\|_{L^\infty} \le C \|\nabla^{s_2} f\|_{L^2}^{\theta} \|\nabla^{s_3} f\|_{L^2}^{1-\theta},$$

*where $0 < \theta < 1$ satisfies*

$$-\frac{s_1}{3} = \left(\frac{1}{2} - \frac{s_2}{3}\right)\theta + \left(\frac{1}{2} - \frac{s_3}{3}\right)(1-\theta).$$

**Lemma A.3** ([18, 20])**.** *Let $s > 0$ and $1 < p < +\infty$. Assume that $q_j, r_j \in (1, +\infty]$ for $j = 1, 2$ such that*

$$\frac{1}{p} = \frac{1}{q_1} + \frac{1}{r_1} = \frac{1}{q_2} + \frac{1}{r_2}.$$

*Then, for $f \in W^{s,q_1} \cap L^{r_2}$ and $g \in W^{s,q_2} \cap L^{r_1}$, we have*

$$\|fg\|_{\dot{W}^{s,p}} \le C(\|f\|_{\dot{W}^{s,q_1}}\|g\|_{L^{r_1}} + \|f\|_{L^{r_2}}\|g\|_{\dot{W}^{s,q_2}}),$$
$$\|fg\|_{W^{s,p}} \le C(\|f\|_{W^{s,q_1}}\|g\|_{L^{r_1}} + \|f\|_{L^{r_2}}\|g\|_{W^{s,q_2}}).$$

*Furthermore, if $f \in \dot{W}^{s,q_1} \cap \dot{W}^{1,r_2}$ and $g \in \dot{W}^{s-1,q_2} \cap L^{r_1}$, we have*

$$\|[f, |\nabla|^s]g\|_{L^p} \le C(\|f\|_{\dot{W}^{s,q_1}}\|g\|_{L^{r_1}} + \|\nabla f\|_{L^{r_2}}\|g\|_{\dot{W}^{s-1,q_2}}).$$

**Lemma A.4.** *For $\ell, \vartheta, s \ge 0$ with $\ell + s + \vartheta > 0$, we have*

$$\|\Lambda^\ell f\|_{L^2} \le \|\Lambda^{\ell+s} f\|_{L^2}^{1-\theta} \|f\|_{\dot{H}^{-\vartheta}}^{\theta}, \text{ where } \theta = \frac{s}{\ell + s + \vartheta}. \tag{A.2}$$

*Proof.* By Parseval's identity, the relation $\ell = (1-\theta)(\ell+s) - \theta\vartheta$, and Hölder's inequality, we have

$$\|\Lambda^\ell f\|_{L^2}^2 = \int_{\mathbb{R}^3} \left(|\xi|^{2(\ell+s)}|\widehat{f}(\xi)|^2\right)^{1-\theta} \left(|\xi|^{-2\vartheta}|\widehat{f}(\xi)|^2\right)^{\theta} \mathrm{d}\xi \le \|\Lambda^{\ell+s} f\|_{L^2}^{2(1-\theta)} \|f\|_{\dot{H}^{-\vartheta}}^{2\theta}.$$

Thus (A.2) follows, and the proof is complete. □

**Lemma A.5** ([32])**.** *Let $0 < \vartheta < 3$ and $1 < p < q < +\infty$ satisfy $1/p = 1/q + \vartheta/3$. Then, we have*

$$\|\Lambda^{-\vartheta} f\|_{L^q} \le C \|f\|_{L^p}.$$

**Lemma A.6** ([27])**.** *Let $\sigma \ge 0$, $1 < p < +\infty$, and let $G : \mathbb{R} \to \mathbb{R}$ be a smooth function satisfying $G(0) = 0$. Then, for any function $u \in L^\infty \cap \dot{W}^{\sigma,p}$, there holds*

$$\|\Lambda^\sigma G(u)\|_{L^p} \le C(\|u\|_{L^\infty}) \|\Lambda^\sigma u\|_{L^p}.$$

**Lemma A.7.** *Let $s_1, s_2 \in \mathbb{R}$ satisfy $s_1, s_2 < 3/2$ and $s_1 + s_2 > 0$. Then for any functions $u \in \dot{H}^{s_1}(\mathbb{R}^3)$ and $v \in \dot{H}^{s_2}(\mathbb{R}^3)$, we have*

$$\|uv\|_{\dot{H}^{s_1+s_2-\frac{3}{2}}} \leq C\|u\|_{\dot{H}^{s_1}}\|v\|_{\dot{H}^{s_2}}. \tag{A.3}$$

*Proof.* The estimate is a standard homogeneous product estimate. It can be derived using Bony's paraproduct decomposition along with the product laws in homogeneous Besov spaces. More precisely, taking $p = 2$ and spatial dimension $N = 3$ in the standard Besov product estimates yields

$$\|uv\|_{\dot{B}^{s_1+s_2-\frac{3}{2}}_{2,2}} \leq C\|u\|_{\dot{B}^{s_1}_{2,2}}\|v\|_{\dot{B}^{s_2}_{2,\infty}}.$$

We refer to [2, 9] for detailed discussions of such Besov product laws. Then (A.3) follows from the standard homogeneous Besov–Sobolev identification and the continuous embedding

$$\dot{H}^s(\mathbb{R}^3) = \dot{B}^s_{2,2}(\mathbb{R}^3), \qquad \dot{B}^s_{2,2}(\mathbb{R}^3) \hookrightarrow \dot{B}^s_{2,\infty}(\mathbb{R}^3), \quad \text{for } s \in \mathbb{R}.$$

However, for the sake of completeness and to keep the argument within the Sobolev framework, we provide a direct proof below using only Sobolev spaces.

Let $s_3 := 3/2 - s_1 - s_2$. By duality, the desired estimate (A.3) is equivalent to the trilinear estimate

$$\left|\int_{\mathbb{R}^3} uvw \,\mathrm{d}x\right| \leq C\|u\|_{\dot{H}^{s_1}}\|v\|_{\dot{H}^{s_2}}\|w\|_{\dot{H}^{s_3}}. \tag{A.4}$$

Notice that the assumptions imply

$$s_1, s_2, s_3 < 3/2, \qquad s_1 + s_2 > 0, \qquad s_1 + s_3 = 3/2 - s_2 > 0, \qquad s_2 + s_3 = 3/2 - s_1 > 0.$$

Consequently, at most one of $s_1, s_2, s_3$ can be negative. We proceed by distinguishing two cases.

**Case 1.** $s_1, s_2, s_3 \geq 0$. Define $p_i \in [2, +\infty)$ by $\frac{1}{p_i} = \frac{1}{2} - \frac{s_i}{3}$ for $i = 1, 2, 3$. Since $\sum_{i=1}^3 s_i = \frac{3}{2}$, we have $\sum_{i=1}^3 \frac{1}{p_i} = 1$. By Hölder's inequality and the Sobolev embedding $\dot{H}^{s_i}(\mathbb{R}^3) \hookrightarrow L^{p_i}(\mathbb{R}^3)$, we obtain

$$\left|\int_{\mathbb{R}^3} uvw \,\mathrm{d}x\right| \leq \|u\|_{L^{p_1}}\|v\|_{L^{p_2}}\|w\|_{L^{p_3}} \leq C\|u\|_{\dot{H}^{s_1}}\|v\|_{\dot{H}^{s_2}}\|w\|_{\dot{H}^{s_3}}.$$

**Case 2.** Exactly one of $s_1, s_2, s_3$ is negative. Without loss of generality, by symmetry we assume $s_3 < 0$. Let $r := -s_3 = s_1 + s_2 - 3/2$. Note that $0 < r < \min\{s_1, s_2\}$ since $s_1, s_2 < 3/2$. By duality,

$$\left|\int_{\mathbb{R}^3} uvw \,\mathrm{d}x\right| \leq \|uv\|_{\dot{H}^r}\|w\|_{\dot{H}^{-r}}.$$

Thus, it suffices to estimate $\|uv\|_{\dot{H}^r}$. Using Lemma A.3, we have

$$\|uv\|_{\dot{H}^r} \leq C\left(\|\Lambda^r u\|_{L^{q_1}}\|v\|_{L^{r_1}} + \|u\|_{L^{r_2}}\|\Lambda^r v\|_{L^{q_2}}\right).$$

The exponents are chosen such that $\frac{1}{q_1} = \frac{s_2}{3}$, $\frac{1}{r_1} = \frac{1}{2} - \frac{s_2}{3}$, $\frac{1}{r_2} = \frac{1}{2} - \frac{s_1}{3}$, and $\frac{1}{q_2} = \frac{s_1}{3}$. This ensures $\frac{1}{q_j} + \frac{1}{r_j} = \frac{1}{2}$ for $j = 1, 2$.

Notice that $s_1 - r = 3/2 - s_2$ and $s_2 - r = 3/2 - s_1$. Applying the Sobolev embeddings $\dot{H}^{s_1-r}(\mathbb{R}^3) \hookrightarrow L^{q_1}(\mathbb{R}^3)$ and $\dot{H}^{s_2}(\mathbb{R}^3) \hookrightarrow L^{r_1}(\mathbb{R}^3)$, we get

$$\|\Lambda^r u\|_{L^{q_1}} \leq C\|u\|_{\dot{H}^{s_1}}, \quad \|v\|_{L^{r_1}} \leq C\|v\|_{\dot{H}^{s_2}}.$$

Similarly, applying the analogous Sobolev embeddings to the second term yields the exact same estimate. Therefore,

$$\|uv\|_{\dot{H}^r} \leq C\|u\|_{\dot{H}^{s_1}}\|v\|_{\dot{H}^{s_2}},$$

which yields (A.4) because $\|w\|_{\dot{H}^{-r}} = \|w\|_{\dot{H}^{s_3}}$.

In both cases, (A.4) holds. The lemma follows directly by duality. ☐

**Lemma A.8.** *Let $r_1, r_2 > 0$ and $r_1 + r_2 > 1$. Then, for any $t \geq 0$, we have*

$$\begin{aligned} &\int_0^{t/2} (1+t-\tau)^{-r_1}(1+\tau)^{-r_2}\,\mathrm{d}\tau + \int_{t/2}^t (1+\tau)^{-(r_1+r_2)}\,\mathrm{d}\tau \\ &\leq C \begin{cases} (1+t)^{-(r_1+r_2-1)}, & 0 < r_2 < 1, \\ (1+t)^{-r_1}\ln(1+t), & r_2 = 1, \\ (1+t)^{-r_1}, & r_2 > 1. \end{cases} \end{aligned} \tag{A.5}$$

*Consequently, for any $0 < r < \min\{r_1, r_1 + r_2 - 1\}$, one has*

$$\int_0^{t/2} (1+t-\tau)^{-r_1}(1+\tau)^{-r_2}\,\mathrm{d}\tau + \int_{t/2}^t (1+\tau)^{-(r_1+r_2)}\,\mathrm{d}\tau \leq C(1+t)^{-r}. \tag{A.6}$$

*If $r_2 \neq 1$, one may also take $r = \min\{r_1, r_1 + r_2 - 1\}$.*

*Proof.* For $0 \leq \tau \leq t/2$, we have $1 + t - \tau \geq (1 + t)/2$, whereas for $t/2 \leq \tau \leq t$, we have $1 + \tau \geq (1 + t)/2$. Therefore,

$$\begin{aligned} &\int_0^{t/2} (1+t-\tau)^{-r_1}(1+\tau)^{-r_2}\,\mathrm{d}\tau + \int_{t/2}^t (1+\tau)^{-(r_1+r_2)}\,\mathrm{d}\tau \\ &\leq 2^{r_1}(1+t)^{-r_1}\left(\int_0^{t/2}(1+\tau)^{-r_2}\,\mathrm{d}\tau + \int_{t/2}^t (1+\tau)^{-r_2}\,\mathrm{d}\tau\right) \\ &= 2^{r_1}(1+t)^{-r_1}\int_0^t (1+\tau)^{-r_2}\,\mathrm{d}\tau. \end{aligned}$$

A direct calculation gives

$$\int_0^t (1+\tau)^{-r_2}\,\mathrm{d}\tau \leq C \begin{cases} (1+t)^{1-r_2}, & 0 < r_2 < 1, \\ \ln(1+t), & r_2 = 1, \\ 1, & r_2 > 1, \end{cases}$$

which proves (A.5). The estimate (A.6) follows by absorbing the logarithmic factor when $r_2 = 1$. □

**Acknowledgments:** J.H.Z. is supported in part by National Science Centre, Poland, grant number 2023/50/A/ST1/00447. He is deeply grateful to Professor Agnieszka Świerczewska-Gwiazda for her kind and generous support during his postdoctoral research. He also sincerely thanks the University of Warsaw for its warm hospitality.

**Conflict of Interest:** The authors declare that they have no conflict of interest.

**Data availability:** Data sharing not applicable to this article as no datasets were generated or analyzed during the current study.

Central China Normal University, School of Mathematics and Statistics, Wuhan, 430079, China
*Email address*: `yuanxu2019@163.com`

University of Warsaw, Institute of Applied Mathematics and Mechanics, Banacha 2, 02-097 Warszawa, Poland
*Email address*: `jhzhang@mimuw.edu.pl`